\documentclass{amsart}
\usepackage{geometry}                
\usepackage{graphicx}
\usepackage{amssymb}
\usepackage{epstopdf,amsmath}
\usepackage{enumitem}

\newtheorem{prop}{Proposition}[section]
\newtheorem{definition}{Definition}[section]

\newtheorem{theorem}[prop]{Theorem}
\newtheorem{conj}[prop]{Conjecture}

\newtheorem*{prop*}{Proposition}
\newtheorem*{theorem*}{Theorem}

\newtheorem*{mainques}{Main Question}

\newcommand{\Poly}{\textrm{Poly}}
\newcommand{\tr}{\textrm{Trace}}

\newcommand{\Prob}{\textrm{Prob}}

\newcommand{\CC}{\mathbb{C}}
\newcommand{\TT}{\mathbb{T}}
\newcommand{\eps}{\epsilon}

\newcommand{\ZZ}{\mathbb{Z}}
\newcommand{\RR}{\mathbb{R}}

\newcommand{\Reb}{\textrm{Re}}
\newcommand{\Imb}{\textrm{Im}}

\title{Large value estimates in number theory, harmonic analysis, and computer science}
\author{Larry Guth}

\begin{document}

\maketitle

\section{Introduction}

 The large value problem for a matrix $M$ asks: if $v$ is an input vector with a given norm,  how many entries of $Mv$ can have size at least some threshold $\lambda$?  Large value problems come up in many parts of math, including analytic number theory, harmonic analysis / PDE, and computer science.   There are many longstanding open questions about them.   In this survey,  we explore what we know about large value problems and what makes them difficult.
 
One of our main examples will be the large value problem for Dirichlet polynomials -- a problem in analytic number theory related to the zeroes of the Riemann zeta function.   This problem boils down to a large value estimate for a particular matrix $M_{Dir}$.   Recently,  James Maynard and I made a little progress on this problem \cite{GM}.   My original plan for this survey was to discuss different methods for the large value problem for Dirichlet polynomials.   

But I gradually decided that it makes sense to focus on large value problems in general.   I believe that most of the issues about the large value problem for Dirichlet polynomials are issues about large value problems in general.   Mathematicians in several different fields have run into these issues and grappled with them.   The issues are by no means resolved,  but I think each field has an interesting perspective on them,  and I will try to share perspectives from analytic number theory,  harmonic analysis, and computer science. 

One of our goals is to survey different methods for studying large value problems.    There are three fundamental techniques to prove large value estimates, and they seem to have been discovered independently by mathematicians in several fields.   We will call them 

\begin{enumerate}

\item  The operator norm method, based on the operator norm of $M$.  

\item The power method.  This is based on taking some type of tensor power of $M$ and studying its operator norm.  

\item The $M M^*$ method.  Based on studying the entries of $M M^*$. 

\end{enumerate}

These techniques are quite general and they give some bounds for the large value problem for many matrices $M$.  Sometimes the bounds are sharp, but often they are not.  In each field, it has been a significant problem to go beyond these techniques,  and in each field people have found a way to do so.  

 In the early 90's,  in \cite{B3},  Bourgain proved new large value estimates for a matrix $M_{S}$ that describes solutions to the Schrodinger equation.  
  This remarkable work is based on a special structure for that matrix,  called a wave packet structure.   This work was followed up by lots more exploration using wave packets in harmonic analysis.   In a striking example from 2014, Bourgain and Demeter \cite{BD} proved {\it sharp} large value estimates for a matrix $M_{PS}$ which describes periodic solutions to the Schrodinger equation.  This work is all based on the special wave packet structure that appears in these matrices,  and the method doesn't adapt to other matrices such as $M_{Dir}$ or matrices coming from computer science.  

Recently, in \cite{GM}, James Maynard and I  improved on the techniques above for Dirichlet polynomials.  At almost the same time, in \cite{DHPT},  Diakonikoas, Hopkins, Pensia, and Tiegel improved on the techniques above for a large value problem from computer science.  Some of the proof ideas are related to each other (and were discovered independently).  One of our main goals is to describe some ideas from these recent works.

Our other main goal is to describe the issues and difficulties in trying to prove large value estimates.  There is a large range of parameters where the operator norm method is the best known method in the large value problem for Dirichlet polynomials.  This proof is based on orthogonality, it takes less than a page, and it was known a hundred years ago.  I find it striking that it is really hard to improve on this one page argument.  The work in \cite{GM} improves on the operator norm method for some new range of the parameters, but there remains today a large range of parameters where this one page argument gives the best known bound.

Given a large matrix $M$,  solving the large value problem for $M$ is a difficult computational problem.   There is no known polynomial time algorithm that gives a good approximation to the answer.   I think that the difficulty of this computational problem is related to the difficulty of proving sharp large value estimates in analytic number theory and harmonic analysis.   So we will survey what is known about this problem from the point of view of computational complexity. 

One interesting direction here concerns a problem about matrices with planted structure,  which we call the planted large value problem.   Let us introduce the problem in the following way.
Analytic number theorists have spent a long time wondering, ``how do we know that there isn't a sneaky vector $v$ so that $M_{Dir} v$ has many large components?  If there were such a sneaky vector $v$,  how would we find it?"   The planted large value problem turns this vague question into a precise question.

Imagine you are given a matrix $M$ that was constructed in one of two ways.   It might be a random matrix $M$,  which obeys very strong large value estimates.  Or it might be a matrix with a planted structure,  designed so that there is a sneaky input vector $v$ such that $Mv$ has many large components.   Since the random matrix does not have any sneaky vector $v$,  you can check whether the given matrix has a planted structure by exhaustive search,  which takes exponential time.   The planted large value problem asks whether we can detect the planted structure in polynomial time.

We will survey what is known about the planted large value problem.   The problem is not fully understood,  but there is an interesting conjectural picture related to the sum of squares hierarchy in computational complexity theory.   Based on ideas in this story,  Mao and Wein \cite{MW} proved that in a certain range of parameters,  the planted large value problem cannot be solved using an algebraic method called low degree testing.
In fact,  they show that in this range of parameters,  low degree testing cannot do better than the simple one page operator norm method discussed above.    There is a conjecture in computer science called the low degree testing conjecture which predicts that for certain types of problems,  no polynomial time algorithm can do better than low degree testing.   So if the low degree testing conjecture is true,  then in this range of parameters, the planted large value problem cannot be solved in polynomial time.  

I think that it is interesting to compare the planted large value problem with the large value problem for Dirichlet polynomials.   We will review the known methods for these two problems and we will see that the methods and the results are quite similar. 

This comparison raises an important question.   The planted large value problem is about comparing a random matrix to a matrix with planted structure.   But matrices that occur in pure math,  such as $M_{Dir}$,  are not random matrices.   They may have special structure,  and that special structure might be key to proving large value estimates.

We will also discuss matrices with special structure and how special structures have been used to prove large value estimates.   In this direction,  we will focus on the work in harmonic analysis using wave packets that we discussed above -- Bourgain's work \cite{B3} on the restriction problem,  and Bourgain-Demeter's work \cite{BD} on decoupling,  which gives {\it sharp} large value estimates for certain special matrices.    This work is 
a striking example of a large value estimate based on the special structure of a particular matrix.

The situation for the matrix $M_{Dir}$ is not so clear.   As we will see,  the matrix $M_{Dir}$ has some special structures that are relevant to the large value problem,  but it's not clear whether these special structures have a major impact on what bounds we can prove.

\vskip10pt

Here is an outline of the essay.

\vskip10pt

In Section 2, we state the large value problem and survey how large value problems appear in different parts of math.

In Section 3,  we work through some examples and use them to motivate the main conjectures about the large value problem for Dirichlet polynomials.

In Section 4,  we give an overview of what is known about approximately solving the large value problem computationally.

Next we discuss the main methods for proving large value estimates.  

In Section 5, we explain the three fundamental techniques for proving large value estimates mentioned above, and we discuss their limitations.

In Sections 6-7, we discuss the recent work in \cite{GM} and \cite{DHPT} that goes a little beyond the techniques from Section 5.  

Next we discuss some barriers to proving sharp large value estimates. 

In Section 8, we discuss large value problems from the point of view of computational complexity,  including the planted large value problem.  

In Section 9,  we discuss another barrier to proving sharp large value estimates for Dirichlet polynomials,  which is related to the Kakeya problem in goemetric measure theory.   This barrier was discovered by Bourgain in \cite{Bour}.

In Section 10,  we discuss some special structures that appear in harmonic analysis and which lead to stronger large value estimates for some special matrices.


\vskip10pt

{\bf Acknowledgements.}  Thanks to James Maynard,  Sam Hopkins,  and Tim Hsieh for interesting conversations about these topics.

\vskip10pt

{\bf Notation. }

\vskip10pt

If $v \in \CC^N$ is a vector,  we write $\| v \|$ for the Euclidean norm

$$ \| v \| = \left( \sum_j |v_j|^2 \right)^{1/2}, $$

and $\| v \|_p$ for the $\ell^p$ norm:

$$ \| v \|_p = \left( \sum_j |v_j|^p \right)^{1/p}.$$

Throughout the paper, we will assume that $M$ is a $T \times N$ matrix and that $N \le T \le N^{O(1)}$.  

We write $A \lesssim B$ to mean that there is a constant $c$ so that $A \le c B$.   

We write $A \lessapprox B$ to mean that for every $\eps > 0$,  there is a constant $c_\eps$ so that $A \le c_\eps N^\eps B$.

\section{Large value problems in different parts of math}

\subsection{Quantitative linear algebra}

The most fundamental object in quantitative linear algebra is the operator norm of a matrix $M$ and the closely related singular value decomposition of the matrix $M$.   Recall that we write $\| v \| = \| v \|_2$ for the standard Euclidean norm of $v$.   The operator norm of $M$ is defined by

\begin{equation}
\label{opernorm} \| M \| =  \max_{\| b \| =1} \| Mb \|. 
\end{equation}

\noindent We write the singular values of $M$ as $s_1(M) \ge s_2(M) \ge ... $ so that $s_1(M) = \| M \|$.  

There is an efficient algorithm to compute $\| M \|$ and to compute the singular value decomposition of $M$.   Partly related to that, we can often estimate the operator norm of matrices that appear in pure math.

The large value problem is part of a circle of more complex questions in quantitative linear algebra.  
Suppose that $M$ is a matrix.   The large value problem asks,  if the input vector $v$ is normalized in some way,  then how many entries of the output $Mv$ can have norm bigger than some threshold.   We make this precise in the following definition.

\begin{definition} \label{LVdef} If $M$ is a $T \times N$ matrix,  then 

$$ LV_{M,  \ell^p}(\lambda) := \max_{\| v \|_p = N^{1/p}} \# \{ i : | (Mv)_i | > \lambda \}. $$
\end{definition}

The normalization $\| v \|_p = N^{1/p}$ is just a convention.   It has the feature that if $|v_j| = 1$ for every $j$,  then $v$ is normalized for every choice of $p$.   

Estimating $LV_{M, \ell^p}(\lambda)$ is more complex than estimating the norm $\| M \|$.  For a large matrix $M$, there is a brute force algorithm to approximate $LV_{M, \ell^p}(\lambda)$ in exponential time, but there is no known polynomial time algorithm that gives a good approximation.  In pure math, there are many specific matrices for which the large value problem is not understood, even after a lot of effort.

Let us briefly explain why large value estimates are more difficult to understand than the operator norm.   Suppose that $W$ is a subset of the rows of $M$,  and suppose we want to know whether there is a vector $b$ with $\| b \|_2$ normalized so that $| (Mb)_j|$ is big for $j \in W$.    Let $M_W$ be the sub-matrix of $M$ that includes only the rows of $W$.    By computing $\| M_W \|$,  we can get good information about the size of $|(Mb)_j|$ for $j \in W$.   However,  we do not know which set $W$ to examine.   To get good information about the large value problem with this method,  we would have to examine every set $W$,  and this takes exponential time.   Optimizing over the choice of $W$ was explicitly considered in \cite{DHPT}.   For any $\eta \in [0,1]$,  they defined the $\eta$-sparse singular value of $M$ as

\begin{equation}
SSV_{\eta}(M) := \max_{|W| \le \eta T} \| M_W \|
\end{equation}

\noindent The sparse singular value problem and the large value problem are very closely related - almost equivalent.  They involve optimizing over the choice of $W$.   There are many choices for $W$ and no obvious shortcut for finding a near optimal choice.   We will come back to this issue at length in the later sections.

A close cousin of the large value problem is estimating the $p$-to-$q$ norm of a matrix.  Recall that the $p$-to-$q$ norm of the matrix $M$ is defined as

$$ \| M \|_{p \rightarrow q} := \max_{\| b \|_{\ell^p} =1} \| Mb \|_{\ell^q}. $$

The norm $\| M \|$ is $\| M \|_{2 \rightarrow 2}$.  But if $p > 1$ and $2 < q < \infty$,  then estimating $\| M \|_{p \rightarrow q}$ is closely related to the large value problem.  In this regime, 
there is no known polynomial time algorithm that gives a good approximation to $\| M \|_{p \rightarrow q}$.  And in pure math, there are many specific matrices for which $\| M \|_{p \rightarrow q}$ is not understood, even after a lot of effort.

In this section, we recall how large value problems come up in several parts of math.  In the rest of the survey, we explore what is known about them.

\subsection{The large value problem for Dirichlet polynomials}

The large problem for Dirichlet polynomials developed as an approach to proving estimates about the zeroes of the Riemann zeta function.   Recall that $\zeta(s)$ is defined by 

\begin{equation} \label{zeta} \zeta(s) = \sum_{n=1}^\infty n^{-s}.\end{equation}

\noindent This sum converges when $\Reb(s) > 1$,  but $\zeta(s)$ extends to a meromorphic function on $\CC$ with a single pole at $s=1$.   The Riemann zeta function vanishes at some points in the critical strip $0 \le \Reb(s) \le 1$.   These points are called the zeroes of the zeta function,  and they encode information about the density of prime numbers.   Suppose we write the zeroes of $\zeta(s)$ as $\rho_j = \sigma_j + i t_j$.    Each zero $\rho_j$ corresponds to a certain part of the fluctuation in the density of primes.   The real part $\sigma_j$ describes the amplitude of the fluctuation and the imaginary part $t_j$ describes the frequency of the fluctuation.   The Riemann hypothesis states that $\sigma_j = 1/2$ for every $j$.   This implies that the fluctuation in the density of primes is small.   The Riemann hypothesis is out of reach,  and so it's interesting to prove weaker bounds on the number of zeroes with large $\sigma$.   Analytic number theorists try to bound $N(\sigma, T)$,  the number of zeroes $\rho_j$ with $\sigma_j \ge \sigma$ and $|t_j| \le T$.   Bounds for $N(\sigma, T)$ lead to estimates about the distribution of primes especially on short intervals.

To try to understand when $\zeta(s) = 0$,  we study $\frac{1}{\zeta(s)}$.   In the range $\Reb(s) > 1$,  $\frac{1}{\zeta(s)}$ can also be described by an infinite series:

\begin{equation} \label{mobius} \frac{1}{\zeta(s)} = \sum_{n=1}^\infty \mu(n) n^{-s},  \end{equation}

\noindent where $\mu(n)$ is the Mobius function.   Recall that $\mu(n)$ depends on the prime factorization of $n$: if $n$ is a product of $r$ distinct primes,  then $\mu(n) = (-1)^r$,  and if $n$ is not a product of distinct primes $\mu(n) = 0$.

This series does not converge if $\Reb(s) \le 1$,  but we can still consider truncated pieces of the form 

\begin{equation} \label{mobiustrunc}
D_N(s) = \sum_{n=N+1}^{2N} \mu(n) n^{-s}.
\end{equation}

\noindent If $\zeta(s) = 0$,  then one of the truncated pieces $D_N(s)$ or a closely related function must be large - see Chapter 10 of \cite{IK} for an introduction to this topic and a precise statement.

The function  (\ref{mobiustrunc}) is an example of a Dirichlet polynomial.   If we want to estimate the number of zeta zeroes with $\sigma_j$ very close to $\sigma$,  then we can suppose $\Reb(s) = \sigma$ and try to understand how (\ref{mobiustrunc}) depends  on $t = \Imb(s)$.    Our goal is to estimate how often this function is large.

The main difficulty is that we can prove very little about the coefficients $\mu(n)$.   Understanding $\mu(n)$ requires understanding how the factorization properties of $n$ are distributed,  which is essentially equivalent to our original problem,  and only very weak bounds are known.   For instance,  it is known that $| \sum_{n=N+1}^{2N} \mu(n) | = o(N)$, but if is far out of reach to prove a bound of the form  $| \sum_{n=N+1}^{2N} \mu(n) | \lesssim N^{.99}$.   Other than very weak bounds,  the main thing we know about $\mu(n)$ is $|\mu(n)| \le 1$ for every $n$.   This has led analytic number theorists to try to prove bounds for a much more general class of Dirichlet polynomials.  

A Dirichlet polynomial of degree $N$ is a trigonometric polynomial of the form

\begin{equation} \label{defdir} D(t) = \sum_{n=N}^{2N} b_n n^{it} \end{equation}

\noindent The large values problem asks about the size of the superlevel set $\{ t: |D(t)| > N^\sigma \}$.  We will formulate the problem precisely as follows.

\begin{mainques} Suppose that

\begin{enumerate}

\item $D(t)$ is a Dirichlet polynomial of degree $N$ as in (\ref{defdir}).

\item $|b_n| \le 1$ for all $n$.

\item $W \subset \ZZ \cap [0,T]$, and $|D(t)| > N^\sigma$ on $W$.  

\end{enumerate}

What is the largest possible cardinality of $W$, in terms of $N$, $T$, and $\sigma$?
\end{mainques}

(The reason that we consider $t \in \ZZ$ is that $|D(t)|$ is morally constant on intervals of length 1, and so understanding $D(t)$ at integer values of $t$ essentially captures the behavior of $D(t)$.)

This question is the large value problem for a certain matrix $M_{Dir}$.   The input vector is $b$ with components $b_n$.   The output vector $M_{Dir} b$ encodes the values of $D(t)$ at $t = 1, ..., T$.   So $M_{Dir}$ is a $T \times N$ matrix and its components are 

$$(M_{Dir})_{t,n} = n^{it}. $$

Note that $(M_{Dir} b)_t = D(t)$.  So the large value problem above is exactly $LV_{M_{Dir}, \infty}(N^\sigma)$.  So to summarize,  large value estimates for the matrix $M_{Dir}$ lead to bounds on the zeroes of zeta.

\subsection{Large value estimates in Fourier analysis and PDE}

We can think of Dirichlet polynomials from the point of view of Fourier analysis as trigonometric polynomials with certain frequencies.   Since $n^{it} = e^{i t \log n}$,  we can rewrite a Dirichlet polynomial $D(t)$ as

$$ D(t) = \sum_{n = N}^{2N} b_n e^{i t \log n}, $$

\noindent a trigonometric polynomial where the frequencies come from the set $\{ \log n \}_{n=N}^{2N}$.  

A number of problems in Fourier analysis fit into this framework.   We can set up a matrix analogous to $M_{Dir}$ to describe the behavior of trigonometric polynomials with any given set of frequencies.  Suppose that $\Phi \subset \RR^d$ is a finite set of frequencies.  A trigonometric polynomial with frequency set $\Phi$ is a sum of the form 

$$u(x) = \sum_{\xi \in \Phi} b_{\xi} e^{ i x \cdot \xi}. $$

Here $x \in \RR^d$.   To get a discrete problem,  we consider a finite set of points $x$.   If we normalize so that $\Phi \subset B_1$,  then it is natural to consider points $x \in \ZZ^d \cap B_R$ for some large radius $R$.   Now the behavior of trigonometric polynomials on the ball of radius $R$ is encoded by a matrix $M_\Phi$.   The columns of $M_\Phi$ correspond to frequencies $\xi \in \Phi$,  and the rows of $M_\Phi$ correspond to points $x \in \ZZ^d \cap B_R$.   The entries of $M_\Phi$ are given by $(M_\Phi)_{x, \xi} = e^{i x \cdot \xi}$,  so that $(M_\Phi b)_x = u(x)$. 

One motivation for problems of this sort comes from PDE.   For any constant coefficient linear PDE on $\RR^d$,  we can write the solutions to the PDE in the form

\begin{equation} \label{fourrep} u(x) =  \int_{S} b(\xi) e^{i \xi \cdot x} d vol_S(\xi). \end{equation}

\noindent Here $S$ is a hypersurface that depends on the given PDE and $d vol_S$ is a smooth measure on $S$.   For example,  we consider the free Schrodinger equation 

$$\partial_t u(y,t) = i \sum_{j=1}^{d-1} \partial_j^2 u(y,t). $$

\noindent Here $y \in \RR^{d-1},  t \in \RR$,  and we set $x = (y,t) \in \RR^d$.   The Schrodinger equation describes a quantum mechanical particle.  
In this case,  the surface $S$ is the paraboloid defined by $\omega_d =  \sum_{j=1}^{d-1} \omega_j^2$.  

The formula (\ref{fourrep}) involves a linear operator $E_S$ defined by

\begin{equation} \label{defEs} E_S b(x)  =  \int_{S} b(\xi) e^{i \xi \cdot x} d vol_S(\xi). \end{equation}

\noindent The linear operator $E_S$ is like an infinite-dimensional matrix,  and it can be approximated by a large matrix $M_{\Phi_S}$ of the type above,  where $\Phi_S$ is a discrete set that approximates the hypersurface $S$.  

Solutions to the Schrodinger equation tend to disperse -- or spread out -- over time.   A solution $u(y,t)$ may be highly concentrated for a short time,  but it cannot remain concentrated for a long time.   Large value estimates for the operator $E_S$ (or the matrix $M_{\Phi_S}$) give quantitative bounds on how much $u(y,t)$ can concentrate in space-time.   Such estimates play an important role in linear and non-linear PDE.  
 The field of restriction theory in Fourier analysis is concerned with large value estimates for $E_S$ and related operators.   We will give a short survey of this field in Section \ref{secharmanal}.

\subsection{Large value estimates in statistics and computer science}

Large value problems appear in computer science in several ways.   They are often related to finding sparse structures in data.

One example is called sparse principal component analysis.  Suppose that $v_n \in \RR^T$ is a sequence of vectors,  $n= 1, ..., N$.   In classical principal analysis,  we look for a unit vector $w$ that has a large correlation with many vectors $v_n$: in other words,  we try to find

\begin{equation} \label{pca} \max_{\| w \| =1 } \sum_{n = 1}^N | v_n \cdot w |^2. 
\end{equation}

\noindent Principal component analysis is a fundamental tool in statistics and in many other areas.   In sparse principal component analysis,  we look for a sparse unit vector $w$ that has a large correlation with many vectors $v_n$.   For $\eta \in [0,1]$,  we say that a vector $w \in \RR^T$ is $\eta$-sparse if the support of $w$ has at most $\eta T$ components.   The sparse PCA problem with sparsity $\eta$ is to try to find

\begin{equation} \label{spca} \max_{\| w \| =1, w \textrm{ is $\eta$-sparse} } \sum_{n = 1}^N | v_j \cdot w |^2. 
\end{equation}

See \cite{XZ} for an overview of statistical applications of sparse PCA.   

Let $M$ denote the $T \times N$ matrix with columns $v_1, .., v_N$.   Principal component analysis is finding the largest singular value and singular vector of the matrix $M$.  There are efficient algorithms to do this.  In contrast, the sparse PCA problem is closely related to the large value problem for $M$, and there is no known polynomial time algorithm that gives a good approximation.

Another example from computer science is trying to find a sparse vector in a subspace.    Suppose we are given an orthonormal basis for a subspace $S \subset \RR^T$ and we would like to decide whether $S$ contains a sparse vector.   Suppose that the dimension of $S$ is $N$ and the orthonormal basis is $v_1, ..., v_N \in \RR^T$.    Again,  we let $M$ be the matrix with columns $v_1, ..., v_N$.   Note that $b_1 v_1 + ... + b_N v_N = Mb$ and that $\sum_{t=1}^T |(Mb)_t|^2 = \sum_{n=1}^N |b_n|^2$.  Therefore,  if $b_1 v_1 + ... + b_N v_N$ is sparse,  then its non-zero entries are large.   So looking for a sparse vector in $S$ is closely related to the large value problem for $M$.

These tasks and other similar tasks occur fairly often in statistics and computer science.   Linear algebra plays a key role in these fields,  and many tasks boil down to quantitative linear algebra.  In many cases, the quantitative linear algebra boils down to the singular value decomposition.  But there are also many cases when it boils down to the large value problem or its cousins -- especially if the problem involves a sparse vector.  The introduction to \cite{DHPT} mentions a few other ways that large value problems appear in computer science.  








\section{Examples and conjectures}

In this section, we start to get intuition about the large value problem by looking at some examples.  With that motivation, we can state the main conjectures about the large value problem for Dirichlet polynomials.

Recall that $M$ is a $T \times N$ matrix and that 

$$LV_{M, \ell^p}(N^\sigma) = \max_{\| b \|_p \le N^{1/p}} \# \{ i: | (Mb)_i | > N^\sigma \}. $$

\noindent So we consider $b$ with $\| b \|_p \le N^{1/p}$,  and we suppose $| (Mb)_i| > N^\sigma$ for $i \in W$.   Then we want to understand $|W|$.

\subsection{Examples}

We begin to study the large value problem by looking at examples.  There are two basic examples that are relevant to the large value problem for any matrix $M$.  

\vskip10pt

{\bf Example 1.}  Random input.  Suppose that $M$ is a $T \times N$ matrix whose entries are unit complex numbers.  Suppose that the input vector $b$ is a random vector where each component $b_n$ is $\pm 1$ with equal probability.  Then for each $t$, $| (Mb)_t |$ is usually $\sim N^{1/2}$.  So if $\sigma = 1/2$, then we can have $|W| \sim T$, which is the maximal possible value.

\vskip10pt
The next example shows that for any $\sigma \le 1$, there are many examples with $|W| \approx N^{2 - 2 \sigma}$.

\vskip10pt

{\bf Example 2.}  Focusing on an arbitrary set.  Suppose that $M$ is a $T \times N$ matrix whose entries are unit complex numbers and $1/2 \le \sigma \le 1$.   Suppose we want to choose a vector $b$ with $|b_n| \le 1$ so that $| (Mb)_t |$ is as large as possible for a given $t$.  The optimal choice is $b_n = \bar M_{tn}$, which gives $(Mb)_t = N$.  In other words, we choose $b$ to be the complex conjugate of the $t^{th}$ row of $M$, which we write as $M_t$.  Now suppose that we want to choose $b$ to make $| (Mb)_t |$ large for $t$ in a set $U$.  A reasonable approach is to choose $b$ to be a linear combination of the vectors $\{ \bar M_{t} \}_{t \in U}$.  This method leads to the following proposition.

\begin{prop} \label{propfocex}
Suppose that $1/2 < \sigma \le 1$.  Suppose that $U \subset \{1, ..., T \}$ is any subset with $|U| \sim N^{2 - 2 \sigma}$.  Then we can choose a vector $b$ so that $|b_n| \lessapprox 1$ and so that $|(Mb)_t| \gtrsim N^\sigma$ for $t \in W$, where $W \subset U$ with $|W| \sim N^{2 - 2 \sigma}$.
\end{prop}

For most matrices $M$, we will see that these two examples are essentially optimal.

\subsection{Random matrices} \label{subsecranmat}

We consider a random $T \times N$ matrix, where each entry is chosen independently from a nice distribution,  such as the uniform distribution on the unit complex numbers,  or the uniform distribution on $\pm 1$,  or the standard Gaussian $N(0,1)$.  
 We can then consider the behavior of a `typical' matrix.  The rough behavior of a typical random matrix is easy to work out, and it serves as a useful reference point when we consider the behavior of a particular matrix such as $M_{Dir}$.  For a typical random matrix, the two examples from the last subsection are essentially optimal.  We make this precise in the following proposition.  

\begin{prop} \label{proprandlv} Suppose that $N \le T \le N^{O(1)}$.   Suppose that $M$ is a random $T \times N$ matrix where each entry is chosen independently from one of the distributions above.    Then with high probability,  $M$ obeys the following large value estimate.

\begin{equation} \label{lvopt}
\textrm{ If $\sum_j |v_j|^2 \le N$,  and if $\sigma > 1/2$,  and if $| \sum_j M_{ij} v_j | \ge N^\sigma$ for $i \in W$,  then $|W| \lessapprox N^{2 - 2 \sigma}$.}
\end{equation}

As a corollary, 

\begin{equation} \label{lvoptlq}
\textrm{ if $\sum_j |v_j|^2 \le N$,  then $\sum_i  |(Mv)_i|^q \lessapprox N^q + N^{q/2} T$.}
 \end{equation}
\end{prop}

(This proposition also applies if each entry is chosen according to another nice distribution with mean zero, such as $\pm 1$ or the Gaussian $N(0,1)$.)

Since most matrices $M$ obey the large value estimate (\ref{lvopt}), it is reasonable to ask whether the matrix $M_{Dir}$ does also.  Montgomery originally conjectured that (\ref{lvopt}) would hold for $M_{Dir}$, but Bourgain found a counterexample to this conjecture.

\subsection{Examples based on approximate geometric series} \label{subsecappgeo}

A geometric series is a trigonometric polynomial of the following form:

$$G(t) =  \sum_{n = 0}^{L-1} e^{ i t (a + d n)} = e^{i a t} \sum_{n = 0}^{L-1} e^{i t d n}. $$

Geometric series play an important role in Fourier analysis in general, and they are important examples for large value problems.  We can sum geometric series exactly, and we get

$$G(t) = e^{i a t} \frac{1 - e^{i t d L}}{1 - e^{i t d}}. $$

\noindent The norm of $G(t)$ is very large when $td$ is close to a multiple of $2 \pi$ and quite small otherwise.   


A Dirichlet polynomial can never be exactly a geometric series, but Bourgain observed that a Dirichlet polynomial can be an approximate geometric series.  This occurs when the coefficients $b_n$ are supported on a small subinterval of $\{ N, ..., 2N \}$ -- say the interval $I$ defined by $I = \{n_I + 1, ..., n_I + L \} \subset [N, 2N]$.  If $L$ is small compared to $N$, then for $n \in I$, $\log n$ is well approximated by its degree 1 Taylor series around $n_I$:

$$ \log (n_I + h) = \log n_I + \frac{1}{n_I} h + O(n_I^{-2} h^2) \approx \log n_I + \frac{1}{n_I} h .$$

\noindent So when $L$ is small, the sequence of frequencies $\{ \log n \}_{n \in I}$ is close to an arithmetic progression, and the Dirichlet polynomial $D(t) = \sum_{n \in I} e^{i t \log n}$ is close to a geometric series.
A normalized version of this example shows that Dirichlet polynomials do not obey the very strong large value estimates (\ref{lvopt}) and (\ref{lvoptlq}),  which hold for random matrices $M$.

\subsection{The main conjectures about Dirichlet polynomials}

In the counterexample based on approximate geometric series, the coefficients $b_n$ are mostly zero.  They are normalized to obey an $\ell^2$ normalization $\sum_n |b_n|^2 = N$, but $\max_n |b_n|$ is much bigger than $1$.   If we insist that $\max_n |b_n| \le 1$,  then we rule out this counterexample. 

In the application to the Riemann zeta function,  we actually need bounds for Dirichlet polynomials obeying the $\ell^\infty$ normalization $\max_n |b_n| \le 1$. 

Following Bourgain's counterexample,  Montgomery made a modified conjecture,  based on the normalization $\max_n |b_n| \le 1$.

\begin{conj}[Montgomery's large value conjecture, cf {M2} page 142] \label{conjMontgomery} Suppose $N \le T \le N^{O(1)}$.  If $\sigma > 1/2$,  then

$$LV_{M_{Dir},  \ell^\infty}(N^\sigma) \lessapprox N^{2 - 2 \sigma}. $$

In other words, if $D(t) = \sum_{n=N}^{2N} b_n e^{i t \log n}$ with $|b_n| \le 1$ and $W \subset \ZZ \cap [0,T] $ such that $|D(t)| > N^\sigma$ for $t\in W$,  then

\[
|W| \lessapprox_\sigma N^{2 - 2 \sigma}.
\]

\end{conj}

This conjecture implies a strong estimate about the zeroes of the Riemann zeta function called the `density hypothesis'.  For some applications, especially involving the distribution of primes in short intervals, the density hypothesis leads to bounds that are comparable to those coming from the full Riemann conjecture.  To prove the density hypothesis, it is enough to prove a weaker bound: $|W| \lessapprox T^{2 - 2 \sigma}$.  This weaker bound is equivalent to a sharp $\ell^q$ estimate for $D(t)$ for all $q$.

\begin{conj}[Montgomery's $\ell^q$ conjecture] \label{conjMontgomerylq} Suppose $N \le T \le N^{O(1)}$.  For every $q \ge 2$, 

$$ \| M_{Dir} \|_{\infty \rightarrow q} \lessapprox N + N^{1/2} T^{1/q}. $$

In other words,  for any $q \ge 2$, if $D(t) = \sum_{n=N}^{2N} b_n e^{i t \log n}$ with $|b_n| \le 1$,  then

$$ \sum_{t =1}^T |D(t)|^q \lessapprox N^q + N^{q/2} T. $$

\end{conj}

In Conjecture \ref{conjMontgomerylq},  for a given $T, N$, there is a critical value of $q$ given by solving $N^q = N^{q/2} T$.   The bound with the critical $q$ implies the bound with all other $q$ and hence the whole conjecture.  

To summarize,  Conjecture \ref{conjMontgomery} implies Conjecture \ref{conjMontgomerylq} which implies the density hypothesis for the zeroes of the Riemann zeta function.  But both Conjecture \ref{conjMontgomery} and Conjecture \ref{conjMontgomerylq} are far out of reach of current techniques.

\subsection{Evidence for the main conjectures}

A few special cases of the main conjectures are known.   We will explain the methods in detail a little later on.

\begin{itemize}

\item If $T = N$,  then both conjectures are known.   This follows from an orthogonality argument, which lets us find the operator norm of the matrix $M_{Dir}$.  

\item The weaker conjecture,  Conjecture \ref{conjMontgomerylq},  is known when $T$ is an integer power of $N$.  In this case, the critical exponent $q$ is an even integer.   The bound then follows by combining a powering trick and orthogonality.

\item  In the late 60s,  Halasz-Montgomery introduced a new method that proves the strong form of the large value conjecture,  Conjecture \ref{conjMontgomery},  when $\sigma$ is sufficiently large.   This work helped motivate Conjecture \ref{conjMontgomery}.  The method is based on estimating the entries of $M_{Dir} M_{Dir}^*$.  

\end{itemize}

Conjectures \ref{conjMontgomery} and \ref{conjMontgomerylq},  in their current form,  date from the late 1980s,  when Bourgain found the example based on approximate geometric series.   In that time,  no one has found a counterexample.   

The other evidence in favor of these conjectures is the large value estimate for random matrices,  Proposition \ref{proprandlv}.   According to this Proposition,  if $M$ is a random $T \times N$ matrix with unit complex entries,  then with high probability,  $M$ obeys the analogue of Conjecture \ref{conjMontgomery}.   This suggests that when we encounter an explicit matrix $M$,  such as $M_{Dir}$,  then we should expect it to obey Conjecture \ref{conjMontgomery} unless we see a good reason,  based on the structure of $M$,  why it should not do so.

\section{Large value estimates from the computational point of view}

In this section, we start to consider the problem of approximating $LV_{M, \ell^p}(\lambda)$ from the computational point of view.   This problem is natural in its own right as a problem in computer science, but I think it also gives some perspective on the difficulty of proving sharp large value estimates in pure math, such as the main conjectures about large value of Dirichlet polynomials, Conjectures \ref{conjMontgomery} and Conjecture \ref{conjMontgomerylq}.  

Recall that if $M$ is a $T \times N$ matrix, then

$$ LV_{M,  \ell^p}(\lambda) := \max_{\| v \|_p = N^{1/p}} \# \{ j : | (Mv)_j | > \lambda \}. $$

Given an arbitrary large matrix $M$,  one can try to approximately solve the large value problem computationally.    If $M$ is a $T \times N$ matrix,  there is a simple brute force algorithm that gives a good approximation to the large value problem in exponential time.    But we do not know any polynomial time algorithm that gives a good approximation of $LV_{M,  \ell^p}(\lambda)$.   It would be interesting to know how accurately a polynomial time algorithm can approximate $LV_{M, \ell^p}(\lambda)$.   The known polynomial time algorithms give upper and lower bounds that differ by a power of $N$.  

It is a closely related problem to approximate $\| M \|_{p \rightarrow q}$ in polynomial time.   In the regime $2 < q < \infty$,  $p > 1$,  the known polynomial time algorithms give upper and lower bounds that differ by a power of $N$ (cf. \cite{GMU} and references therein).   For some values of $p,q$ in this range,  it is known to be $NP$-hard to approximate $\| M \|_{p \rightarrow q}$ within a constant factor -- see \cite{BV},  \cite{BBHKSZ},  and \cite{BGGLT}.    The methods in these papers may lead to similar results for the large value problem.

Because of these computational difficulties,  
there is no meaningful numerical evidence supporting the main conjectures about large values of Dirichlet polynomials.   This is different from the situation with the Riemann hypothesis itself.  There is a lot of numerical evidence supporting the Riemann hypothesis -- for instance, we know that the first $10^{9}$ zeroes of the Riemann zeta function have real part equal to 1/2 (EDIT: Reference).  In contrast, we are not able to check the large value conjectures numerically even for $N = 200$.

There is not currently a conjectural picture of how well polynomial time algorithms can approximate $LV_{M, \ell^p}(\lambda)$.  If $M$ is allowed to vary among all $T \times T$ matrices,  then it looks reasonable to conjecture that it is $NP$-hard to approximate the function $LV_{M, \ell^p}(N^\sigma)$ to within a factor $T^{\gamma}$, for some $\gamma  > 0$.  But there is no clear conjecture about the sharp value of the exponent $\gamma$.

Proving bounds for $LV_{M, \ell^p}(\lambda)$ is not just hard in the worst case, it is also hard in the average case.  If $M$ is a random $T \times N$ matrix, then Proposition \ref{proprandlv} shows that with high probability, $M$ obeys an essentially optimal large value estimate.  But if we are given a particular matrix $M$ that was sampled from the random distribution in Proposition \ref{proprandlv}, it is still hard to prove that $M$ obeys strong large value estimates.  In computer science, this is called the problem of certifying large value estimates for $M$.  It is closely related to the problem of distinguishing an honest random matrix $M$ from a matrix that has a planted structure which causes $LV_{M, \ell^p}(\lambda)$ to be much larger.

I think the problem of certifying large value estimates for random matrices is closely parallel to the large value problem for Dirichlet polynomials.  The two problems have been studied independently, but the methods and bounds that are known for the two problems are closely parallel, as we will see in the next few sections.

The problem of certifying large value estimates for random matrices is part of the field of average case computational complexity.   For a broad class of related problems in the field, the best current method is called the sum of squares hierarchy.  The limits of this method are well understood in a number of cases.  It is not fully understood what large value estimates for a random matrix $M$ can be certified by the sum of squares method, but there was a lot of recent progress in \cite{DHPT}.  Based on recent developments in the field, there is good evidence that the sum of squares method cannot certify that a random matrix obeys the large value estimates in Proposition \ref{proprandlv} or Conjecture \ref{conjMontgomery} or \ref{conjMontgomerylq}.  

Experts in the field consider it to be plausible that the sum of squares method is the best polynomial time algorithm for a broad class of problems in average case computational complexity, including certifying large value estimates for random matrices.  So there is a plausible conjectural picture about what large value estimates can be certified in polynomial time.  In this picture, almost every random matrix obeys very strong large value estimates, but these strong large value estimates cannot be certified by a polynomial time algorithm. 

After this brief discussion of computational complexity, let us return to the question why the large value problem for Dirichlet polynomials is hard. Recall that $M_{Dir}$ is the $T \times N$ matrix that encodes Dirichlet polynomials.   We would like to prove bounds for $LV_{M_{Dir}, \ell^\infty}(\lambda)$ for all $N$ and all $T = N^{O(1)}$.  However, let us set a more modest goal and try to prove bounds when $N= 10^3$ and $T = 10^4$.  In this case, $M_{Dir}$ fits quite easily in a computer.  Let us say that we would like our proof to have length at most $10^{25}$, so that there is some hope of checking it on a large computer.  Even this problem looks very hard, and there is no plausible approach to it in the literature on large value problems.  If we replace $M_{Dir}$ by a random matrix $10^4 \times 10^3$ matrix $M$, then there is no evidence that any such proof exists. 

This discussion raises the question whether $M_{Dir}$ behaves like a random matrix.  Perhaps $M_{Dir}$ has some special structure that makes large value bounds easier to prove than for a random matrix.  Before we discuss this point, let us consider the problem of designing a matrix $M$ that can be proven to obey large value estimates.

Conjecture \ref{conjMontgomery} holds for a random matrix $M$ with high probability, and so it holds for many matrices $M$ with entries of norm 1.  However, no specific matrix $M$ with entries of norm 1 has been proven to obey Conjecture \ref{conjMontgomery}.  

This problem is closely related to the restricted isometry problem from compressed sensing.  If $M$ is a matrix and $W$ is a subset of the rows of $M$, recall that $M_W$ is the submatrix consisting of those rows.  Let us say that $M_W$ is an almost-isometry if the singular values of $M_W$ all agree up to a factor of $1.01$.   A matrix $M$ has the restricted isometry property up to rank $S$ if $M_W$ is an almost isometry whenever $ |W| \le S$.   Suppose that $T = N^\alpha$ for some $1 < \alpha < 2$.   A random $ T \times N$ matrix is known to have the restricted isometry for all $S \le N^{1 - \epsilon}$,  But no specific matrix $M$ with entries of norm 1 has been proven to have the restricted isometry property up to this rank.  The restricted isometry property is closely related to the large value problem, and the methods in the literature are related.  For more information about the restricted isometry problem and the best current bounds for a specific matrix, see \cite{BDFKK}.

Conjecture \ref{conjMontgomerylq} is weaker than Conjecture \ref{conjMontgomery},  and it has been proven for some specific matrices.   Recall that Conjecture \ref{conjMontgomerylq} is about estimating $\| M_{Dir} \|_{\infty \rightarrow q}$ for all $q \ge 2$.   For any given $T, N$,  there is a critical value of $q$,  $q_{crit}$,  determined by $N^{q_{crit}} = N^{q_{crit}/2} T$.   The difficulty of this problem depends on whether $q_{crit}$ is an even integer.  When $q_{crit}$ is an even integer,  then we will see a short proof of Conjecture \ref{conjMontgomerylq}.   There are similar methods for many other matrices,  including random matrices.

When $q_{crit}$ is not an even integer,  then it is much more difficult to prove the analogue of Conjecture \ref{conjMontgomerylq}.   Only in the last decade,  such a bound was proven for a specific family of matrices with entries of norm 1.   The result was proven by Bourgain and Demeter \cite{BD} in their breakthrough work on decoupling.   The matrix they studied encodes periodic solutions to the Schrodinger equation,  and the bound proven by Bourgain and Demeter is called the (sharp) Strichartz estimate for the periodic Schrodinger equation.   We will discuss their work later in Section \ref{secharmanal}, and we will see that their proof is based on special structural properties of this particular matrix.   

Now we come back to the question whether $M_{Dir}$ has special structure that makes large value estimates for $M_{Dir}$ easier to prove than for a random matrix.  The matrix $M_{Dir}$ does have some special structures that have been used to study the large value estimate, and we will review these structures during the survey.  However, the impact of these special structures so far has not been that large, and the general shape of bounds that we can currently prove for $M_{Dir}$ is similar to the bounds we can prove for a random matrix in polynomial time.  The periodic Schrodinger equation studied by Bourgain-Demeter has very special structure which has a major impact on understanding large value estimates and $p$-to-$q$ norms for that matrix.  The matrix $M_{Dir}$ also has some special structures that have been studied, but it is not clear whether these special structures have a major impact on understanding large value estimates.  

In the next sections, we discuss the known methods to prove large value estimates.  After that, in Section \ref{secbarrierlowdeg}, we return to the computational complexity of the large value problem and flesh out this discussion with more details.

\section{Classical methods for studying large value problems} \label{secclassmeth}

In this section,  we will study three fundamental methods of studying the large value problem for a matrix $M$.  These methods have been discovered independently in several communities in different contexts.  

\begin{enumerate}

\item The operator norm method.

\item The power method.

\item The $M M^*$ method.

\end{enumerate}

Before we get started, let's recall the problem.  

\begin{mainques} (For Dirichlet polynomials) Suppose that

\begin{enumerate}

\item $D(t) = \sum_{n=N+1}^{2N} b_n e^{i t \log n}$.

\item $|b_n| \le 1$ for all $n$.

\item $W \subset \ZZ \cap [0,T]$, and $|D(t)| > N^\sigma$ on $W$.  

\end{enumerate}

What is the largest possible cardinality of $W$, in terms of $N$, $T$, and $\sigma$?
\end{mainques}

We can ask a similar problem about any matrix.  

\begin{mainques} (For general matrices) Let $M$ be a $T \times N$ matrix.   Suppose that

\begin{enumerate}

\item $b \in \RR^N$ with $\| b_n \|_{\ell^p} \le N^{1/p}$.  

\item $W \subset \ZZ \cap [1,T]$, and $|(Mb)_t| > N^\sigma$ on $W$.  

\end{enumerate}

For a given matrix $M$,  what is the largest possible cardinality of $W$, in terms of $p$ and $\sigma$?
\end{mainques}

The main question about Dirichlet polynomials is a special case of the second question where $M = M_{Dir}$ is defined by $(M_{Dir})_{t,n} = e^{i t \log n}$ and $p = \infty$.  

We will focus here on proving large value estimates for $M_{Dir}$ or for a random matrix $M_{ran}$.

\subsection{Orthogonality / operator norm method}

Given a matrix $M$, one fundamental thing that we can compute in polynomial time is its operator norm

\begin{equation}
\| M \| = \| M \|_{2 \rightarrow 2 } := \max_{\| v \|_{\ell^2} = 1 } \| M v \|_{\ell^2}.
\end{equation}

For the matrix encoding Dirichlet polynomials, the operator norm can be approximated by a simple orthogonality argument, giving the following $\ell^2$ bound.

\begin{prop} \label{proporthintro} If $D(t) = \sum_{n=N}^{2 N} b_n e^{i t \log n}$ and $T \ge N$, then

$$ \sum_{t = 1}^T |D(t)|^2 \lesssim T \sum_n |b_n|^2. $$

\end{prop}

\begin{proof} [Proof idea] For each $n$, think of $e^{i t \log n}$ as a function on $\{1, ..., T \}$.  These functions are approximately orthogonal to each other.  If they were exactly orthogonal, then we would have $ \sum_{t = 1}^T |D(t)|^2 = T \sum_n |b_n|^2. $
\end{proof}

Similarly,  if $M$ is a random $T \times N$ matrix with entries of norm 1,  then $\| M \| \sim T^{1/2}$ with high probability,  and so we have

$$ \sum_{t = 1}^T | (Mb)_t |^2 \lesssim T \sum_n |b_n|^2. $$

Under the hypotheses in our main question,  we have $\sum_n |b_n|^2 \le N$ and $|D(t)| > N^\sigma$ for $t \in W \subset \{ 1, ..., T \}$),  and so Proposition \ref{proporthintro} gives the bound

\begin{equation} \label{basicorthintro}
|W| N^{2 \sigma} \lesssim T N.
\end{equation}

This basic bound holds for $M_{Dir}$ and for random $M$ and can be checked for many other explicit $M$. 

This bound shows that our main conjectures hold in the regime $T = N$.

\subsection{Power method}

The product of two Dirichlet polynomials is again a Dirichlet polynomial.   To see this,  note that $\log n_1 + \log n_2 =\log (n_1 n_2)$, and then expand out:

$$\left( \sum_{n_1} a_{n_1} e^{i t \log n_1} \right) \left( \sum_{n_2} b_{n_2} e^{i t \log n_2} \right) = \sum_{n} \left( \sum_{n_1 n_2 = n} a_{n_1} b_{n_2} \right) e^{i t \log n}. $$

The final expression is itself a Dirichlet polynomial,  $\sum_n c_n e^{i t \log n}$,  where

$$ c_n = \sum_{n_1 n_2 = n} a_{n_1} b_{n_2}. $$

In particular,  if $D$ is a Dirichlet polynomial of degree $N$ and $k$ is a natural number,  then $D^k$ is a Dirichlet polynomial of degree $N^k$.   Applying Lemma \ref{proporthintro} to $D^k$ leads to sharp bounds for  $\sum_{t=1}^T |D(t)|^{2k}$,  which proves Conjecture \ref{conjMontgomerylq} when $q$ is an even integer.   In fact,  it gives an even better bound where the right-hand side involves $\sum_n |b_n|^2$ instead of $\max_n |b_n|$.  

\begin{prop} \label{propstronglq}
 If $q$ is an even integer,  and if $\sum_n |b_n|^2 \le N$, then

\begin{equation} \label{stronglq}
\sum_{t=1}^T \left| \sum_{n \sim N} b_n e^{it \log n} \right|^q \lessapprox N^q + N^{q/2} T.
\end{equation}

\end{prop}

Note that this method gives sharp bounds when the exponent $q$ is an even integer but it does not apply for any other exponent $q$.   In many places in Fourier analysis,  even integer exponents are easier to understand than other exponents because of this powering trick.

Barak,  Brandao,  Harrow, Kelner, Steurer and Zhou developed a version of the power method for random matrices in \cite{BBHKSZ}.  
For random matrices $M$ it also gives sharp bounds for $\sum_i |(Mv)_i |^{2k}$.    We illustrate the method for $k=2$.   We know that $(Mv)_i = \sum_j M_{ij} v_j$.   Therefore,  we have

\begin{equation} \label{matrixsq} \left( (Mv)_i \right)^2 = \left( \sum_j M_{ij} v_j \right)^2 = \sum_{j_1, j_2} M_{i j_1} M_{i j_2} v_{j_1} v_{j_2}. 
\end{equation}

We can rewrite this equation in terms of a new matrix $M^{\otimes 2}$,  which is a kind of tensor square of $M$.   The rows of $M^{\otimes 2}$ are indexed by $i$ and the columns are indexed by pairs $(j_1, j_2)$.   We have

$$ M^{\otimes 2}_{i,  (j_1, j_2)} = M_{i j_1} M_{i j_2}. $$

Similarly we define 

$$ v^{\otimes 2}_{(j_1, j_2)} = v_{j_1} v_{j_2}. $$

Then we can rewrite equation (\ref{matrixsq}) as

$$  \left( (Mv)_i \right)^2 = \left( M^{\otimes 2} v^{\otimes 2} \right)_i. $$

Now we square both sides and sum over $i$ to get

$$ \sum_i \left| (Mv)_i \right|^4 = \| M^{\otimes 2} v^{\otimes 2} \|_{\ell^2}^2. $$

We can bound the right-hand side using the operator norm of $M^{\otimes 2}$ and noting that $\| v^{\otimes 2} \|_{\ell^2}^2 = \| v \|_{\ell^2}^4$.  We get

$$ \sum_i \left| (Mv)_i \right|^4 \le  \| M^{\otimes 2} \|_{2 \rightarrow 2}^2  \| v \|_{\ell^2}^4. $$

For any matrix $M$,  the matrix norm $\| M^{\otimes 2} \|$ can be approximated in polynomial time.   For the matrix $M_{Dir}$,  this estimate is equivalent to the discussion above.   The fact that the product of two Dirichlet polynomials is a Dirichlet polynomial gives a nice clean way of estimating $ \| M_{Dir}^{\otimes 2} \|$.   For a random complex matrix $M_{ran}$,  $ \| M_{ran}^{\otimes 2} \|$ obeys the same bounds,  and this gives a $\Poly(N)$ length proof that $M_{ran}$ obeys (\ref{lvoptlq}),  which matches Proposition \ref{propstronglq}.

However,  \cite{BBHKSZ} observed that for a random real matrix $M_{ran}$,  $\| M_{ran}^{otimes 2} \|$ is too large to give sharp bounds for $\| M_{ran} \|_{2 \rightarrow 4}$.   Recall that $M^{\otimes 2}_{i, (j_1, j_2)} = M_{i, j_1} M_{i, j_2}$.    When $j_1 = j_2$,  we have $ M_{i,(j,j)} = M_{i,j}^2 \ge 0$.   If $M_{i,j} = \pm 1$ randomly,  then we have $M_{i, (j,j)} = 1$ for all $i,j$.   Let $w_{j_1, j_2}$ be the vector defined by

$$ w_{j_1, j_2} = \begin{cases} 1 & j_1 = j_2 \\  0 & j_1 \not= j_2 \ \end{cases}. $$

\noindent Then $\| M^{\otimes 2} w \|$ is quite large, and $w$ is nearly the largest singular vector of $M^{\otimes 2}$.   In this case \cite{BBHKSZ} break up the matrix $M^{\otimes 2}$ as a sum of a simple part that captures this large singular vector and a leftover part with small norm.   This method gives a polynomial length proof that a random real matrix $M_{ran}$ obeys the estimate in Proposition \ref{propstronglq},  giving an essentially sharp bound for $\| M_{ran} \|_{2 \rightarrow q}$ when $q$ is an even integer.

\subsection{The Montgomery-Halasz / $M M^*$ method}

 In the late 60s,  Montgomery \cite{M}, building on ideas of Halasz and Turan (cf. \cite{H} and \cite{HT}), proved strong large value estimates for Dirichlet polynomials  when $\sigma$ is large.  For sufficiently large values of $\sigma$,  he proved sharp bounds, matching Conjecture \ref{conjMontgomery}.    (These bounds helped to motivate the large value conjecture,  Conjecture \ref{conjMontgomery}.)
 
 The method applies to any matrix $M$.  For a random $T \times N$ matrix $M$, the Halasz-Montgomery method gives a sharp large value estimate for $\sigma > 3/4$ (matching the bound in Proposition \ref{proprandlv}).  The proof has length $Poly(N)$.  
 
 Similarly,  if we want to prove Conjecture \ref{conjMontgomery} for a fixed $N$ and $T = N^{O(1)}$ and we allow a proof of length $Poly(N)$,  then (assuming a standard conjecture in analytic number theory),  the Halasz-Montgomery method would give a $Poly(N)$ length proof for a sharp large value estimate for $\sigma > 3/4$.  
 
But for $\sigma \le 3/4$, the Halasz-Montgomery method does not give any information about the large value problem at all.

Closely related to methods that have been applied to a number of problems in harmonic analysis and in computer science.  As far as I know, the method was first used by Kolmogorov-Seliverstov \cite{KS} and Plessner \cite{P} in their work on pointwise convergence of Fourier series.  

The method is based on studying the entries of $M M^*$, so we call it the $M M^*$ method.  (In harmonic analysis, it is usually called the $T T^*$ method, because the role of the matrix $M$ is played by a linear operator $T$.)

Here is the basic idea of the method.   Suppose that $M$ is a $T \times N$ matrix and that $W$ is a subset of the rows.   Let $M_W$ be the submatrix of $M$ given by including the rows of $W$.    Understanding the large value problem is closely related to estimating

$$ \max_{|W| = S} \| M_W \| $$

\noindent for each choice of $S$.   (Recall that $\| M_W \| = \| M_W \|_{2 \rightarrow 2}$ is the operator norm of $M_W$.)

We study this problem by looking at $M_W M_W^*$.   First note that $\| M_W M_W^* \| = \| M_W \|^2$.   Also $M_W M_W^*$ is closely related to $M M^*$.   Note that $M M^*$ is a $T \times T$ matrix, and $M_W M_W^*$ is the submatrix where the row and column belong to $W$.   

In the Montgomery-Halasz method,  we estimate the entries of $M M^*$ and use them to bound $\| M_W M_W^* \|$.    Recall that 

$$(M M^*)_{i, i'} = \sum_j M_{i,j} \bar M_{i',j}. $$

\noindent If $M$ has entries of norm 1,  then each diagonal entry of $M M^*$ will be equal to $N$.   The off-diagonal entries typically have smaller norm because of cancellation.   If $M$ has random entries of norm 1,  then the off-diagonal entries of $M M^*$ will have norm $\lessapprox N^{1/2}$.  Conjecturally,  this also holds for $M_{Dir}$.   Using these bounds on the entries of $M_W M_W^*$ leads to sharp estimates for the large value problem when $\sigma > 3/4$ but no information when $\sigma \le 3/4$.   Let's state the relevant estimate as a Proposition.

\begin{prop} Suppose that $M$ is a $T \times N$ matrix where each entry has norm 1,  and suppose that each off-diagonal entry of $M M^*$ has norm $\lessapprox N^{1/2}$.   Suppose that $\sum_n |b_n|^2 \le N$ and that $|(Mb)_t| \ge N^\sigma$ for $t \in W$.   If $\sigma > 3/4$,  then

$$ |W| \lessapprox N^{2 - 2 \sigma}. $$

\end{prop}

\begin{proof}  We break up $M_W M_W^*$ as a diagonal matrix plus an off-diagonal matrix.   Because each entry of $M$ has norm 1,  the diagonal part of $M_W M_W^*$ is just $N I$, and so it has norm $N$.   The off-diagonal matrix is a $|W| \times |W|$ matrix where the entries have norm $\lessapprox N^{1/2}$,  so it has norm $\lessapprox |W| N^{1/2}$.   All together we have

\begin{equation} \label{MWnorm}
\| M_W \|^2 = \| M_W M_W^* \| \lessapprox N + |W| N^{1/2}
\end{equation}

Therefore we have

$$ \sum_{t \in W} | (Mb)_t |^2 \le \| M_W \|^2 \|b \|_{\ell^2}^2 \lessapprox ( N + |W| N^{1/2} ) (\sum_n |b_n|^2 ). $$

If we plug in that $ | (Mb)_t | > N^\sigma$ for $t \in W$ and $\sum_n |b_n|^2 \le N$,  then we get

$$ |W| N^{2 \sigma} \lessapprox N^2 + |W| N^{3/2}. $$

If the first term dominates,  we have $|W| \lessapprox N^{2 - 2 \sigma}$.   If the second term dominates,  we have $N^{2 \sigma} \lessapprox N^{3/2}$ and so $\sigma \le 3/4$.  

\end{proof}

For the matrix $M_{Dir}$,  it is conjectured that all off-diagonal entries have size $\lessapprox N^{1/2}$,  but this is a deep open problem,  closely related to the Lindelof hypothesis.   The known bounds for the off-diagonal entries are weaker but still useful: using van der Corput estimates, Montgomery checked that the off-diagonal terms have norm $\lessapprox T^{1/2}$.   Montgomery's large value estimate uses this weaker bound.   It still gives a sharp large value estimate if $\sigma$ is sufficiently large,  but not for the whole range $\sigma > 3/4$.   Since then,  there have been several further improvements in the range $\sigma > 3/4$.   

This $M M^*$ method has been applied to many problems in harmonic analysis.  Later in this survey,  in Section \ref{secharmanal},
we will discuss a few of them.   In some of these problems,  there is a regime corresponding to $\sigma > 3/4$ where the method gives sharp estimates and a regime corresponding to $\sigma \le 3/4$ where the method gives no information.    For these problems,  we will describe what is known in the regime corresponding to $\sigma \le 3/4$.   Sometimes there is a phase transition in behavior when $\sigma$ crosses 3/4 and sometimes there is not.   Sometimes we know the optimal answer in all regimes but often we do not.

\subsection{On the regime $\sigma \le 3/4$}

In the regime $\sigma \le 3/4$,  we have so far seen only two methods to study the large value problem: the operator norm method and the power method.  When  
 $\sigma \le 3/4$ and  $N < T \le N^{3/2}$,  the power method does not help and the best estimate so far comes the basic operator norm bound  (\ref{basicorthintro}).
 
In this regime, until recently, the basic bound (\ref{basicorthintro}) was the best known bound for the large value problem for Dirichlet polynomials.  In \cite{GM}, Maynard and I improved on the bound (\ref{basicorthintro}) in part of this range.

Similarly for random matrices $M$, in this regime, until recently, the basic bound (\ref{basicorthintro}) was the best known bound for the large value problem that could be certified in polynomial time.  In \cite{DHPT}, 
 Diakonikolas, Hopkins, Pensia, and Tiegel gave a polynomial time algorithm that certifies a bound stronger than (\ref{basicorthintro}) in part of this range.
 
  I find it striking that it is difficult to improve on the basic bound (\ref{basicorthintro}).  The proof of  (\ref{basicorthintro}) is less than a page, and it is about a hundred years old.  
 But when $\sigma < 3/4$ and $T \le N^{3/2}$ it was open for a long time to improve on this bound.  And in a slightly smaller regime, like $\sigma < 7/10$ and $T \le N^{6/5}$, it remains an open problem to improve on this one page argument.

 One reason for the difficulty is an almost counterexample.  There is a cousin problem where the bound (\ref{basicorthintro}) is actually sharp when $\sigma \le 3/4$.   This example helps to understand why it is hard to improve on previous bounds when $\sigma \le 3/4$.   

\subsection{An almost counterexample} \label{subsecac}

In this section,  we present a cousin of the large value problem for Dirichlet polynomials.   The cousin problem shares many features of the original problem,  but for the cousin problem,  the simple orthogonality bound (\ref{basicorthintro}) is sharp when $\sigma \le 3/4$.   This is the first of several `barriers' we will discuss in the survey.   It shows that certain methods cannot make progress on the large value problem.

Yuqiu Fu, Dominique Maldague, and I found this cousin problem in \cite{FGM},  in the context of trying to make progress on the Dirichlet large value problem using ideas related to decoupling.  The cousin problem was later simplified by Alex Cohen.  Decoupling theory is a powerful tool for giving $L^p$ estimates and large value estimates in some questions in Fourier analysis.  However, the methods from decoupling do not distinguish the almost counterexample from the Dirichlet large value problem, and therefore they cannot give any new bounds for $\sigma \le 3/4$.

\newtheorem*{setupac}{Cousin of the main question}

\begin{setupac} Suppose that

\begin{enumerate}

\item $\tilde D(t) = \sum_{n=N}^{2N} b_n e^{i t \sqrt{\frac{n}{N}}}$.  

\item $\sum_{n=N}^{2N} |b_n|^2 \le N$.

\item $W \subset [0,T]$ is a 1-separated set where $|D(t)| \gtrsim N^\sigma$.

\end{enumerate}

What is the largest possible cardinality of $W$, in terms of $N$, $T$, and $\sigma$?

\end{setupac}

The orthogonality estimate in Proposition \ref{proporthintro} generalizes to this setup, and so we get the bound (\ref{basicorthintro}): $N^{2 \sigma} |W| \lesssim T N$.   If $\sigma \in [1/2, 3/4]$, then this bound is actually sharp.

\begin{prop} \label{propac} If $1/2 \le \sigma \le 3/4$, then there is a function $\tilde D(t)$ and set $W$ obeying the hypotheses above with $N^{2 \sigma} |W| \sim T N$.
\end{prop}

\begin{proof} We illustrate with $\sigma = 3/4$.  The key point is that the set of frequencies, $\{ \sqrt{ \frac{n}{N}} \}_{n=N}^{2N}$, contains an arithmetic progression of length $\sim N^{1/2}$.  This arithmetic progression comes from the case where $n = m^2$.  If $n = m^2$, then $\sqrt{ \frac{n}{N}} = \frac{m}{\sqrt{N}}$.  Here $m$ can be any integer in the range $[\sqrt{N}, \sqrt{2N}]$, and so we get an arithmetic progression of length $\sim N^{1/2}$.  Then we define the coefficients $b_n$ by

$$ b_n = \begin{cases}  N^{1/4} & \textrm{ if } n = m^2 \\ 0 & \textrm{ else } \\                        \end{cases}. $$

\noindent With these coefficients, $\tilde D(t)$ becomes a geometric series:

$$ \tilde D(t) = N^{1/4} \sum_{m \in [\sqrt{N}, \sqrt{2N}]} e^{i t \frac{m}{\sqrt{N}}}. $$

Note that $\tilde D(0) \sim N^{3/4}$ and that $\tilde D$ is $2 \pi \sqrt{N}$-periodic.  We define $W = 2 \pi \sqrt{N} \ZZ \cap [0, T]$.  

If $\sigma \le 3/4$, we can make a similar construction using a shorter arithmetic progression.

\end{proof}

Remark: The set $W$ can be taken to be a set of integers: the integers closest to $2 \pi \sqrt{N} \ZZ \cap [0,T]$.  

This almost counterexample implies that some proof methods cannot improve on the bound (\ref{basicorthintro}) when $\sigma \le 3/4$.   For instance,  the Halasz-Montgomery method does not distinguish the almost counterexample from the problem of actual Dirichlet polynomials.  Since this is an important point,  let us explain it a little more carefully.   The Dirichlet polynomial problem is associated with the $T \times N$ matrix $M_{Dir}$,  and it is conjectured that the off-diagonal entries of $M_{Dir} M_{Dir}^*$ have norm $\lessapprox N^{1/2}$.   Similarly,  the cousin problem is associated to a $T \times N$ matrix $M_{AC}$.   I think it is natural to conjecture that the off-diagonal entries of $M_{AC} M_{AC}^*$ also have norm $\lessapprox N^{1/2}$.   If these conjectures are true,  then the Halasz-Montgomery method will give the same results for both problems.   In both cases,  it will give a $\Poly(N)$ length proof for the sharp large value estimate for $\sigma > 3/4$ and no information for $\sigma \le 3/4$.    For the matrix $M_{AC}$,  there would be a phase transition when $\sigma$ crosses 3/4.   For $\sigma > 3/4$,  the sharp bound would be $|W| \lessapprox N^{2 - 2 \sigma}$,  while for $\sigma \le 3/4$,  the sharp bound is $|W| \lessapprox T^2 N^{-2 \sigma}$.   








In the next sections, we will explain some of the ideas from \cite{GM} and \cite{DHPT} that lead to some further progress on large value problems in this regime.

\section{The Schatten tensor} \label{secschatten}

In \cite{GM}, Maynard and I made progress on the large value problem for Dirichlet polynomials for $\sigma$ a little less than 3/4.   In \cite{DHPT},  Diakonikoas, Hopkins, Pensia, and Tiegel made progress on certifying large value estimates for random matrices for $\sigma$ a little less than 3/4.  Some of the proof ideas are related to each other (and were discovered independently).

Let us first state the results a little more precisely.  The main theorem in \cite{GM} applies to the case $T = N^{6/5}$.  

\begin{theorem} \label{maingm} (Theorem 1.1 in \cite{GM}) Recall that $M_{Dir}$ is the $T \times N$ matrix from Section 1.  Suppose that $T = N^{6/5}$ and $\sigma \le 3/4$.

Suppose that $|b_n| \le 1$ for all $n$ and that $| (M_{Dir} b)_t | \ge N^\sigma$ for all $t \in W \subset \{ 1, .., T \}$.  Then

$$ |W| \lessapprox N^{\frac{18}{5} - 4 \sigma}. $$

\end{theorem}

This theorem improves on the basic bound (\ref{basicorthintro}) when $\sigma > 7/10$.  If $\sigma = 3/4$, then the basic bound (\ref{basicorthintro}) gives $|W| \lessapprox N^{7/10}$, the bound in \cite{GM} gives $|W| \lessapprox N^{6/10}$, and Conjecture \ref{conjMontgomery} predicts that $|W| \lessapprox N^{5/10}$.  These bounds for the case $T = N^{6/5} $ also imply new bounds for some other values of $T$.

The results in \cite{DHPT} work equally well for all values of $T$ from $N$ to $N^2$.  If you plug in $T = N^{6/5}$, the bounds in \cite{DHPT} match the bounds in \cite{GM}.

\begin{theorem} \label{maindhpt} (Theorem 4.4 in \cite{DHPT}) Suppose that $N \le T \le N^2$ and suppose that $M$ is a random $T \times N$ matrix, where each entry is chosen independently from the standard Gaussian $N(0,1)$.

For any $\epsilon > 0$, there is a polynomial time algorithm for proving bounds for the large value problem for $M$ with running time $N^{D(\epsilon)}$.  With high probability, this algorithm will prove that $M$ obeys the following bounds.

If $\| b \|_2^2 \le N$, and if $| (Mb)_t | \ge N^\sigma$ for all $t \in W$, and if $\sigma \le 3/4$, then 

$$ |W| \le C(\epsilon) N^\epsilon N^{3 - 4 \sigma} T^{1/2}. $$

\end{theorem}

As mentioned above, if we plug in $T = N^{6/5}$, then the bound in Theorem \ref{maindhpt} from \cite{DHPT} matches the bound in Theorem \ref{maingm} from \cite{GM}.  For other values of $T$, the bounds in \cite{DHPT} are stronger than the bounds in \cite{GM}.  If $T = N^\alpha$, then the bound in Theorem \ref{maindhpt} improves on the basic bound (\ref{basicorthintro}) when $\sigma > 1 - \alpha / 4$.  

The Schatten tensor is a key ingredient in both proofs, discovered independently by the two groups.

Recall that $M$ is a matrix.  If $W$ is a subset of the rows of $M$, then we write $M_W$ for the minor of $M$ containing the rows in $W$.  Our goal is to understand the norm of $M_W$ as a function of $W$.  
One classical way to find the norm of a matrix $M$ is to look at $(M^* M)^r$ or $(M M^*)^r$ for a large power $r$.   We can make this relationship precise using singular values.   
We let $s_1(M) \ge s_2(M) \ge ... $ denote the singular values of $M$.    Recall that the norm of $M$ is $\| M \| = s_1(M)$, and recall that for $p \ge 1$, the Schatten $p$-norm of $M$ is

\begin{equation} \label{defSchattp} \| M \|_{S^p} := \left( \sum_i s_i(M)^p \right)^{1/p}
\end{equation}

Finally recall that for any natural number $r$, we have

$$ \tr \left[ (M^*M)^r \right] = \sum_i s_i(M)^{2r} = \| M \|_{S^{2r}}^{2r}. $$

If we are given a matrix $M$ and a fixed $W$,  then we can efficiently compute $\tr \left[ (M_W^* M_W)^r \right]$ and use it to estimate $\| M_W \|$.    As $r$ increases,  this gives us more and more accurate estimates for $\| M_W \|$.  The main problem is that we want to estimate the maximum of $\| M_W \|$ over all $W$ of a given size $S$.
Unless $S$ is very small,  there are too many possibilities to check them one at a time in total time $\Poly(N)$.  
So we have to understand how $ \tr \left[ (M_W^* M_W)^r \right] $ depends on $W$.  

When we expand out the expression $ \tr \left[ (M_W^* M_W)^r \right] $, it naturally has the form of an $r$-linear tensor.  Recall that an $r$-linear tensor on a complex vector space $V$ is a map $A: V^r \rightarrow \CC$ that is linear in each component.  In our case, the vector space $V$ is $\CC^T$ and we associate each subset $W \subset \{1, ..., T \}$ with the vector $1_W$ whose components are $(1_W)_i = 1$ if $i \in W$ and $(1_W)_i = 0$ if $i \notin W$.   If we just expand out $  \tr \left[ (M_W^* M_W)^r \right] $,  then we find that there is an $r$-linear tensor $S_{M, r}$ so that

$$  \tr \left[ (M_W^* M_W)^r \right]  = S_{M, r}(1_W, ..., 1_W). $$

We call the tensor $S_{M, r}$ the Schatten $r$-tensor of the matrix $M$.

The paper \cite{GM} carefully studies the tensor $S_{M_{Dir}, 3}$.  The paper \cite{DHPT} carefully studies the tensor $S_{M_{ran}, r}$ for all $r \ge 2$. 

For $r=2$, $S_{M, 2}$ is a 2-linear tensor, which means it is a matrix.  To study this matrix, one might find its operator norm or its singular value decomposition.  However, studying the matrix $S_{M, 2}$ is closely related to the $M M^*$ method, and it doesn't lead to any new estimates for the large value problem.  So we have to study the tensor $S_{M, r}$ for $r \ge 3$.

There is a natural notion of the norm of an $r$-linear tensor for any $r$:

\begin{equation} \label{defnormtensor} \| S \| := \max_{\| v_1 \|_{\ell^2} = ... = \| v_r \|_{\ell^2} = 1} | S(v_1, .., v_r) |. \end{equation}

In studying a tensor, it is usually a good first step to estimate its norm.   But in some cases, we can get a better understanding by breaking up the tensor into pieces.  We will write $S_{M,r}$ as a sum of a simple tensor and a leftover piece, and then bound the norm of the leftover piece.   Let us present one such decomposition, following \cite{GM}.  To motivate the decomposition, let us try to make an educated guess about the value of $\sum_i s_i(M_W)^{2r}$ for a typical set $W$ of a given size.  Recall that the Hilbert-Schmidt norm of a matrix $M$ is given by

$$ \| M \|_{HS}^2 = \sum_i s_i(M)^2 = \sum_{i,j} |M_{i,j}|^2. $$

\noindent So if every entry of $M$ has norm 1, then

$$ \sum_i s_i(M_W)^2 = |W| N. $$

\noindent The rank of $M_W$ is at most $|W|$.  If $|W| \le N$, then it might be reasonable to guess that $M_W$ has $|W|$ singular values each of norm around $\sqrt{N}$.  So we would guess that $\sum_i s_i(M_W)^{2r} \approx |W| N^r$.  Indeed, Holder's inequality shows that we always have $\sum_i s_i(M_W)^{2r} \ge |W| N^r$, and in some regimes this is pretty close to an equality.

With this in mind, we can break up the tensor $S_{M,r}$ as follows.  We let $I$ be the $r$-linear tensor with components 

$$ I_{i_1, ..., i_r} = \begin{cases} 1 & \textrm{ if } i_1 = ... = i_r \\ 0 & \textrm{ else }
\end{cases} $$

\noindent It is easy to check that $I(1_W, ..., 1_W) = |W|$.  Now we break up $S_{M, r}$ as

\begin{equation} S_{M,r} = N^r I + S^\Delta_{M, r}
\end{equation}

The paper \cite{GM} uses this decomposition.  We have $\| N^r I \| = N^r$, and \cite{GM} proves that in some regimes, $\| S^\Delta_{M_{Dir}, r}\|$ is much smaller.  

The paper \cite{DHPT} uses a more refined decomposition of the tensor $S^{M, r}$, breaking off several other pieces which are all simpler than $S_{M,r}$ itself.

The next step is to estimate the norm of the leftover piece $S^\Delta_{M,r}$ (or the related tensors in the refined decomposition from \cite{DHPT}).   It's important to note that when $r \ge 3$, estimating the norm of a tensor is much more difficult than estimating the norm of a matrix.  For instance, when $r \ge 3$, there is no known efficient algorithm for approximating the norm of an $r$-tensor.

One standard method of getting a bound for the norm of an $r$-linear tensor is to forget some of the structure and think of the $r$-linear tensor as a 2-linear tensor.  This method is called `flattening the tensor'.  We can think of an $r$-linear tensor $S$ as a linear map from the tensor power $V^{\otimes r} \rightarrow \CC$.  If $r = r_1 + r_2$, then we can think of $S$ as a map $\tilde S: V^{\otimes r_1} \times V^{\otimes r_2} \rightarrow \CC$.   We can now think of $\tilde S$ as a matrix, where the rows are indexed by a basis of $V^{\otimes r_1}$ and the columns are indexed by a basis of $V^{\otimes r_2}$.   Because $\tilde S$ is a matrix,  we can efficiently compute its operator norm $\| \tilde S \|$.   Moreover we have $\| S \| \le \| \tilde S \|$ because $\| \tilde S \| = \max_{t_i \in V^{\otimes r_i}, \| t_i \|_{\ell^2} = 1}|  \tilde S(t_1, t_2)|$ whereas $\| S \|$ is given by taking the maximum over the smaller set where each $t_i$ is a tensor product of unit vectors.  

(While this method is simple, it is actually the best known method for estimating the operator norm of tensors in many situations.  For instance, if $r$ is an even integer, and $S$ is a random $r$-tensor where each entry is independently chosen from a standard Gaussian $N(0,1)$, then this method is essentially the best known polynomial length method of upper bounding $\| S \|$.  We will discuss this more in Section \ref{secbarrierlowdeg}.)

The polynomial time algorithm in \cite{DHPT} computes the norms of the flattened matrices corresponding to the tensors in their decomposition of $S_{M,r}$.  When $M$ is chosen randomly, the paper \cite{DHPT} analyzes the norm of each of these matrices, allowing them to analyze what bounds their polynomial time algorithm will give.  (The degree of the polynomial time algorithm depends on $r$, and as $\epsilon \rightarrow 0$ in Theorem \ref{maindhpt}, they take $r \rightarrow \infty$.)  Since $M$ is random, these auxiliary matrices are also random, but their entries are not independent.  This makes it challenging to analyze the norm of such an auxiliary matrix.  
 The introduction to \cite{DHPT} has a nice overview of the ideas of the proof.

Now let us return to $M_{Dir}$.  Suppose for a moment that we fix $T$ and $N$ and we are looking for a $\Poly(N)$ length proof of a large value estimate for $M_{Dir}$.  The method we have just described using the Schatten tensor gives a proof of length $\Poly(N)$ of some large value estimates for $M_{Dir}$.  It is not known what bounds this method gives for the matrix $M_{Dir}$.  My best guess is that if we applied the method from \cite{DHPT} to $M_{Dir}$, the bounds that it would give are similar to the bounds that it gives for a random matrix.  However, it looks difficult to prove this rigorously.  For some perspective, let's go back to the $M M^*$ method, where we want to estimate the off diagonal entries of $M M^*$.  For a random $T \times N$ matrix $M$, all these off-diagonal entries have norm $\lessapprox N^{1/2}$.  It is conjectured that for $M_{Dir}$, the off-diagonal entries of $M_{Dir} M_{Dir}^*$ also have norm $\lessapprox N^{1/2}$, but this old conjecture is out of reach.  In an analogous way, I would guess that the norms of the auxiliary matrices associated to $M_{Dir}$ behave similarly to the norms of the auxiliary matrices associated to $M_{ran}$, but this naive conjecture is also out of reach.


One issue is the almost counterexample we discussed in Section \ref{subsecac}.  Let $M_{AC}$ be the matrix from the almost counterexample.  Everything we wrote in this section so far applies equally to $M_{Dir}$ and $M_{AC}$.  In order to prove a bound for $M_{Dir}$ which is stronger than the one for $M_{AC}$, we could try to prove a bound for $\| \tilde S_{M_{Dir}, r} \|$ which is much less than $\| \tilde S_{M_{AC}, r} \|$ for some $r$.  In time $\Poly(N)$ we can compute these two norms to high precision, and then I believe we would find that  $\| \tilde S_{M_{Dir}, r} \|$ is much less than $\| \tilde S_{M_{AC}, r} \|$.    But for a proof,  we still need to put our finger on a property distinguishing $M_{Dir}$ from $M_{AC}$ that could be used to bound these norms.  

Eventually, Maynard and I did find a special property of $M_{Dir}$ that distinguishes it from $M_{AC}$ and we were able to use this to prove Theorem \ref{maingm}.  We explain these ideas in the next section.

\section{Bounds using special features of Dirichlet polynomials}

In this section, we outline the main ideas of the proof of Theorem \ref{maingm} from \cite{GM}.

One of the key issues is to distinguish the large value problem for Dirichlet polynomials from the almost counterexample in Section \ref{subsecac}. 
If we compare the setup of our main question with the setup of the almost counterexample, there are two differences.

\begin{enumerate}

\item The frequencies in the trigonometric polynomial $D(t)$ are $\{ \log n \}_{n= N}^{2N}$ instead of $ \{ \sqrt{ \frac{n}{N}} \}_{n=N}^{2N}$.

\item The coefficients obey $|b_n| \le 1$ for all $n$ instead of just $\sum_{n=N}^{2N} |b_n|^2 \le N$ . 

\end{enumerate}

In order to improve on the simple bound (\ref{basicorthintro}) for our main question when $\sigma \le 3/4$, the proof needs to use at least one of these two features.  The proof in \cite{GM} uses both features, although either one by itself may be enough.

Suppose that $|b_n| \le 1$ for all $n$ and that $| (M_{Dir} b)_t) | \ge N^\sigma$ for $t \in W \subset \ZZ \cap [0, T]$.  We want to bound $|W|$.  

There are two cases depending on whether the set $W$ has a lot of additive structure.  Notice that in the almost counterexample, the set $W$ was an arithmetic progression, which has a lot of additive structure.  Sets with additive structure have a special relationship with Fourier analysis, and they often come up as examples for the large value problem for a matrix $M$ that comes from Fourier analysis.  We measure additive structure using the additive energy.  For a finite set $W \subset \ZZ$,  we define the additive energy as

\begin{equation} \label{defE}
E(W) := \# \{ t_1, t_2, t_3, t_4 \in W: t_1 + t_2 = t_3 -t_4 \}
\end{equation}

\noindent The additive energy lies between $|W|^2$ and $|W|^3$.   Sets with lots of additive structure,  such as arithmetic progressions,  have $E(W) \approx |W|^3$.   On the other hand,  random sets typically have $E(W) \approx |W|^2$.   (A bit more precisely,  if $W \subset \{1,  ..., T\}$ is a random set of cardinality at most $T^{1/2}$,  then with high probability $E(W) \lesssim |W|^2$.)   The additive energy is closely related to Fourier analysis.  For a set $W \subset \ZZ$, define

\begin{equation} \label{defhatW} \hat W(\xi) = \sum_{t \in W} e^{i \xi t}
\end{equation}

The additive energy can be rewritten in terms of $\hat W$:

\begin{equation} \label{energyfourier} E(W) = \frac{1}{2 \pi} \int_0^{2 \pi} |\hat W(\xi)|^4 d \xi
\end{equation}

If the set $W$ has small additive energy, then we will use the Schatten tensor method to estimate $|W|$.   A key observation is that the Schatten tensor $S_{M_{Dir}, r}$ has a special structure: if we rewrite the Schatten tensor $S_{M_{Dir}, r}$ using the Fourier transform, then it is surprisingly sparse.   For comparison, the tensor $S_{M_{AC}, r}$ is not sparse.   This part of the proof distinguishes $M_{Dir}$ from $M_{AC}$.  In this part of the proof, we only require $\sum_n |b_n|^2 \le N$, not $|b_n| \le 1$ for all $n$.  

If the set $W$ has large additive energy, then we use work of Heath-Brown from the late 70s, \cite{HB}.  This part of the proof uses that $|b_n| \le 1$ for every $n$.  

For the purposes of understanding the main ideas, we focus on two special cases.

\begin{itemize}

\item The case with no additive structure.  In this case $E(W) \approx |W|^2$.  Sometimes we consider an even more special case when $|\hat W(\xi)| \approx |W|^{1/2}$ except when $|\xi| \lessapprox 1/T$.  

\item The case with maximal additive structure.  In this case $E(W) \approx |W|^3$.  

\end{itemize}

 We will explain each of these cases.  There is actually one further idea in \cite{GM}, needed to handle intermediate energies, but we leave it out of this survey.

\subsection{The case with no additive structure} \label{subsecnoaddstruc}

In this section, we use the Schatten tensor approach,  and we find a special feature of the matrix $M_{Dir}$ and the tensor $S^\Delta_{M_{Dir}, 3}$.  This special structure depends on the fact that the frequencies in a Dirichlet polynomial are exactly $\log n$, and it distinguishes honest Dirichlet polynomials from the almost counterexample.

Suppose $M = M_\Phi$ is the matrix encoding trigonometric polynomials with frequency set $\Phi$.  So the rows of $M_{\Phi}$ are indexed by $t \in \{1, ..., T \}$ and the columns of $M_{\Phi}$ are indexed by $\xi \in \Phi$, and the entry $(M_\Phi)_{t, \xi} = e^{i t \xi}$.  In this situation, $\tr \left[ (M_{\Phi, W}^* M_{\Phi, W})^r \right]$ has a special form which can be described in terms of the Fourier transforms of the set $W$.   We write

$$ \hat W(\xi) = \sum_{t \in W} e^{- i t \xi}.$$

We know that $\sum_i s_i( M_{\Phi, W})^{2r} = \tr \left[ (M_{\Phi, W}^* M_{\Phi, W})^r \right]$.  Writing down the formula for $M_{\Phi, W}$ and expanding out the formula for the trace, we get

\begin{equation} \label{momsingr}
\sum_i s_i( M_{\Phi, W})^{2r} = \sum_{\xi_1, ..., \xi_r \in \Phi} \hat W(\xi_1 - \xi_2) \hat W(\xi_2 - \xi_3) ... \hat W(\xi_r - \xi_1).
\end{equation}


This formula describes the tensor $S_{M_\Phi, r}(1_W, 1_W, ..., 1_W)$, but we are writing the tensor in terms of $\hat W$ instead of $W$.  
When $r=2$,   the formula (\ref{momsingr}) simplifies to 

$$\sum_{\xi_1, \xi_2 \in \Phi} |\hat W(\xi_1 - \xi_2)|^2, $$

\noindent and the terms in the sum are all positive.   This sum measures the size of $\hat W$ on the difference set $\Phi - \Phi$.  But when $r \ge 3$,  the terms of (\ref{momsingr}) are not all positive and it is necessary to show that there is significant cancellation between them in order to get new large value estimates in the range $\sigma \le 3/4$.    This is tricky because we need to prove a bound for all sets $W$.   The case $r=2$ involves the difference set $\Phi - \Phi$,  and the case of larger $r$ involves a generalization which we call a cyclic difference set of order $r$:

\begin{equation} \label{defcycdiff} D^r \Phi := \{ ( \xi_1 - xi_2,  \xi_2 - \xi_3, ..., \xi_r - \xi_1): \xi_1, ..., \xi_r \in \Phi \}. 
\end{equation}

Then we can rewrite (\ref{momsingr}) as

\begin{equation} \label{momsingr2}
\sum_i s_i( M_{\Phi, W})^{2r} = \sum_{(\tau_1, ..., \tau_r) \in D^r \Phi} \hat W(\tau_1) \hat W(\tau_2) ... \hat W(\tau_r) = \sum_{\tau \in D^r \Phi} \hat W^{\otimes r} (\tau).
\end{equation}

(In this setup,  $D^r \Phi$ is really a multi-set -- a set with multiplicities -- and the sum above includes the multiplicities.)

The way the set $D^r \Phi$ is distributed plays a role in bounding this sum.   When $\Phi = \{ \log n \}_{n=N}^{2N}$,  the set $D^r \Phi$ is distributed in a special way which will allow us to prove cancellation in the formula above.   We can see this special structure already for the difference set $\Phi - \Phi$.   To get a perspective,  let us compare honest Dirichlet polynomials with the almost counterexample.  We write

$$ \Phi_{Dir} := \{ \log n \}_{n = N}^{2N}, \Phi_{AC} := \left\{ \sqrt{\frac{n}{N}} \right\}_{n=N}^{2N}. $$

\noindent Let us compare the difference sets $\Phi_{AC} - \Phi_{AC}$ and $\Phi_{Dir} - \Phi_{Dir}$.  

Since $| \Phi_{AC} | = N$,  $| \Phi_{AC} - \Phi_{AC} | = N^2$.   In the difference set $\Phi_{AC} - \Phi_{AC}$,  the number 0 appears with multiplicity $N$,  and there are no other elements in a neighborhood around 0 of radius $\sim 1/N$.   Other than that,  I believe that the set $\Phi_{AC} - \Phi_{AC}$ appears quite random.   For instance,  the typical spacing between consecutive elements is $\sim 1/ N^2$,  but by the birthday paradox we expect to see some very nearby elements with spacing $\sim 1/N^4$.   

For our purposes,  it's important to understand the density of the set $\Phi_{AC} - \Phi_{AC}$ at the scale $1/T$,  which we will write as $\rho_{1/T,  AC}$.   You can think of this density as smoothing out the set $\Phi_{AC} - \Phi_{AC}$ at scale $1/T$.   Technically, this is defined as follows.   Let $\eta_{1/T}$ be a smooth function which is roughly constant on $B_{1/T}$ and rapidly decacying off of it, and normalized to have integral 1.   Then we have

$$\rho_{1/T,  AC}(\xi) := \sum_{\tau \in \Phi_{AC} - \Phi_{AC}} \eta_{1/T}(\xi - \tau). $$

\noindent We will be interested in $T$ significantly bigger than $N$ but much smaller than $N^2$,  say $T = N^{6/5}$.  
The density $\rho_{1/T, AC}$ has a large spike near 0.   Other than that,  it roughly follows a smooth profile but with a significant amount of random noise.   Here is a rough sketch.

\includegraphics[scale=.8]{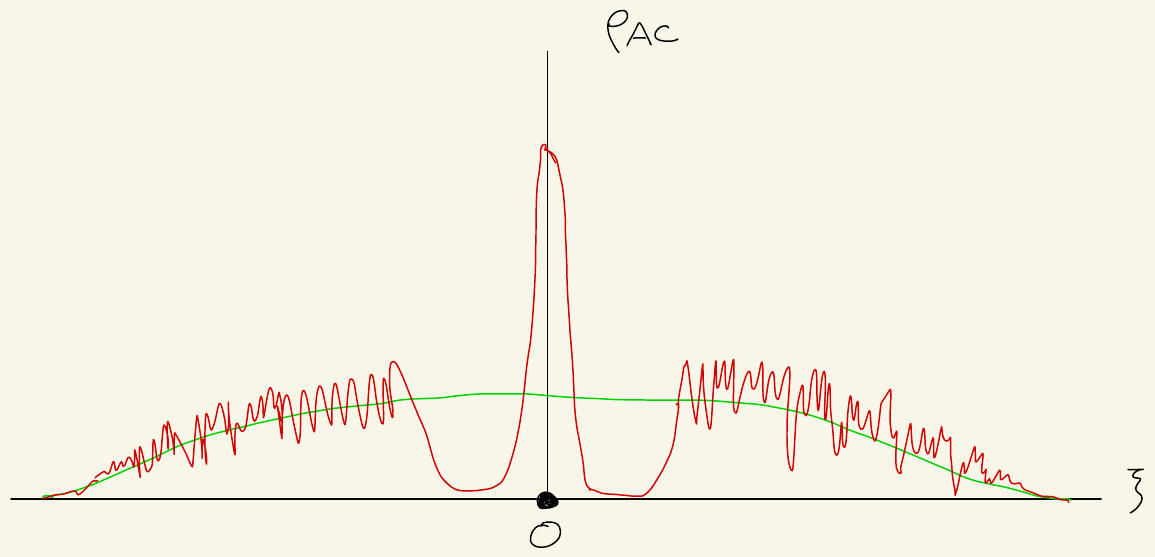}

\newcommand{\smooth}{\textrm{smooth}}
\newcommand{\spikes}{\textrm{spikes}}
\newcommand{\noise}{\textrm{noise}}
\newcommand{\fluc}{\textrm{fluc}}

We can write

\begin{equation} \label{smoothnoise}  \rho_{1/T,  AC} = \rho_{\smooth} + \rho_{\fluc}
\end{equation}

We have $\rho_{\smooth} \sim N^2$.   At 0,  $\rho$ is much larger,  $\sim NT$, and so $\rho_{\fluc}(0) \sim NT$.   At other places,  I believe that $\rho_{\fluc}$ resembles random noise.  For instance,  I expect that (away from 0) $\rho_{\fluc}$ has size $\sim T \sqrt{ N^2/ T} \sim N T^{1/2}$,  which would happen if we replaced $\Phi_{AC} - \Phi_{AC}$ with a random set of $N^2$ points.   This could be checked in time $\Poly(N)$,  but it looks very hard to prove it rigorously.  

Now let's turn to $\Phi_{Dir} - \Phi_{Dir}$.   This is the set of differences $\log n_1 - \log n_2 = \log (n_1/n_2)$.    The spacing of $\Phi_{Dir} - \Phi_{Dir}$ is closely related to the spacing of rational numbers $n_1/n_2$ with $N \le n_1, n_2 \le 2N$.
As above,  the element 0 has multiplicity $N$.   The element $\log (3/2)$ has large multiplicity also because there are many ways to write $3/2 = n_1/n_2$ with $N \le n_1, n_2 \le 2N$.  In general,  we will see large multiplicities at $\log (p/q)$ where $p, q$ are small integers.   But in between these high multiplicity points,  $\Phi_{Dir} - \Phi_{Dir}$ is distributed more evenly than a random set.   For instance,  there is no birthday paradox: the smallest spacing between two distinct elements is comparable to the average spacing,  $\sim 1/N^2$. 

The density $\rho_{1/T,  Dir}$ has multiple spikes corresponding to $\log (p/q)$ where $p,q$ are small integers.   Quantitatively,  the spikes occur when $p,q \lessapprox T/N$.   In between these spikes,  $\rho_{1/T, Dir}$ is extremely smooth.    Here is a rough sketch.

\includegraphics[scale=.8]{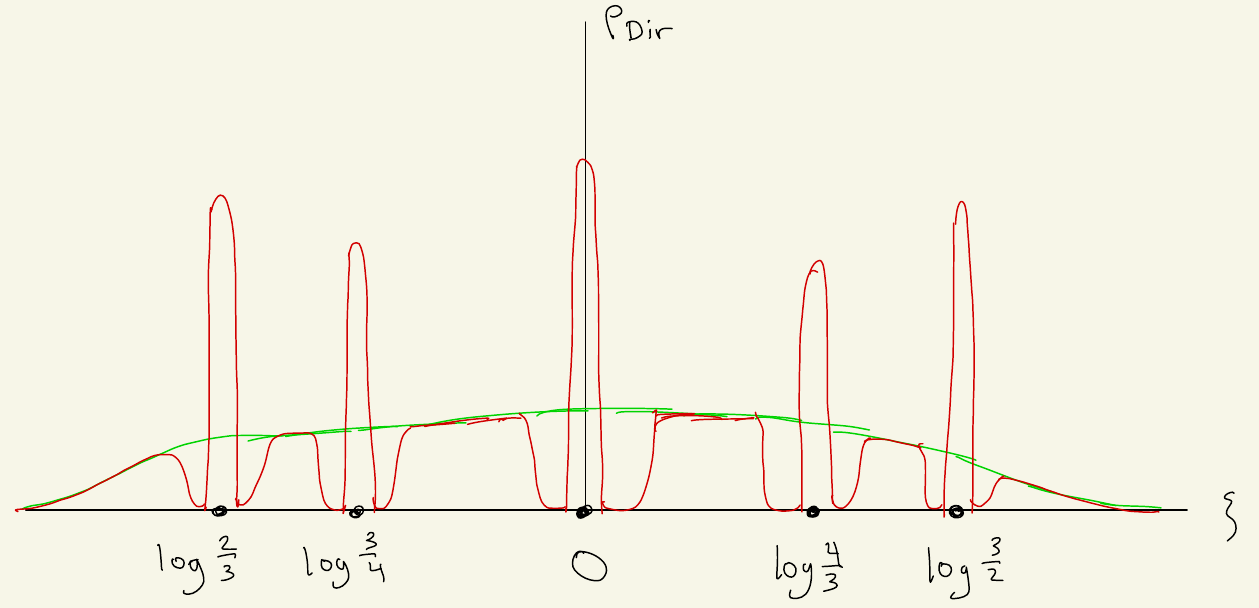}

We can write

\begin{equation} \label{smoothspike}  \rho_{1/T,  Dir} = \rho_{\smooth} + \rho_{\spikes}. 
\end{equation}

Here $\rho_{\smooth}$ is a smooth function of height $\sim N^2$ like above.    The function $\rho_{\spikes}$ is quite large near to $\log (p/q)$,  where $p,q$ are small integers,  but it is essentially zero away from those points.   The spikes occur at $\log p/q$ where $|p|, |q| \lessapprox T/N$ and each spike has width $1/T$.   So if $T$ is, for example,  $N^{6/5}$,  we see a fairly sparse set of tall spikes.  The height of a typical spike is $\sim N^2$,  although the spikes corresponding to small $p,q$ are taller.  

The sets $D^r \Phi_{Dir}$ for larger $r$ behave in a similar way.   For instance,  if $r=3$,  then $D^3 \Phi$ is a subset of a 2-dimensional disk of radius $\sim 1$.   If we take its density at scale $1/T$, then we get a smooth part plus spikes along a small number of curves.   Each curve corresponds to a rational point of the form $(\frac{p_1}{p_3},  \frac{p_2}{p_3})$,  with $|p_i| \lessapprox T/N$,  and so the number of curves is $\lessapprox (T/N)^3$.   The spike over a typical curve has height $\sim N^2$ over a $1/T$-neighborhood of the curve.

These plots of densities are closely related to estimating the singular values of $M_{\Phi, W}$.  Recall that we have (\ref{momsingr})

\begin{equation*}
\sum_i s_i( M_{\Phi, W})^{2r} = \sum_{(\tau_1, ..., \tau_r) \in D^r \Phi} \hat W(\tau_1) \hat W(\tau_2) ... \hat W(\tau_r) 
\end{equation*}

Since $W$ is contained in $[0,T]$, $\hat W$ is morally roughly constant at scale $1/T$, and so this sum is approximately an integral with the density at scale $1/T$: $\sum_i s_i(M_{Dir, W})^{2r}$ is approximately

\begin{equation} \label{momsingrint}
\int  \hat W(\tau_1) \hat W(\tau_2) ... \hat W(\tau_r) \rho_{1/T, Dir}, \end{equation}

which we can break up as

\begin{equation}  \int  \hat W(\tau_1) \hat W(\tau_2) ... \hat W(\tau_r) \rho_{\smooth} + \int  \hat W(\tau_1) \hat W(\tau_2) ... \hat W(\tau_r) \rho_{\spikes}  \end{equation}

Because $\rho_{\smooth}$ is very nice and explicit, we can approximate the first term very accurately.  We do this using Plancherel.  The Fourier transform of $\hat W$ is just $W$, and the Fourier transform of $\rho_{\smooth}$ is highly localized and explicit.  This term works out to $N^r |W|$.   Recall from Section \ref{secschatten} that we broke up $S_{M, r}= N^r I + S^\Delta_{M, r}$.  The contribution of $\rho_{smooth}$ corresponds to the $N^r I$ term, and the contribution of $\rho_{spikes}$ corresponds to the leftover piece $S^\Delta_{M, r}$.  

We bound the contribution of $\rho_{\spikes}$ using the triangle inequality:

$$ \int  \hat W(\tau_1) \hat W(\tau_2) ... \hat W(\tau_r) \rho_{\spikes} \le \int |\hat W(\tau_1) ... \hat W(\tau_r)| |\rho_{spikes}|. $$

\noindent We are able to get a better bound for $M_{Dir, W}$ than for $M_{AC, W}$ because $\int |\rho_{spikes}|$ is much smaller than $\int |\rho_{noise}|$.  If we are assuming that $|\hat W(\tau)| \lessapprox |W|^{1/2}$ (except for $\tau$ close to zero), then we can bound the contribution of $\rho_{spikes}$ by $\lessapprox |W|^{r/2} \int |\rho_{spikes}|$.  


In general, we get good bounds for the second term unless $\hat W(\tau_1) ... \hat W(\tau_r)$ is highly concentrated on the spikes of $\rho_{spikes}$.  In this latter situation, it must happen that $\int_{|\tau| \lesssim 1}  |\hat W(\tau)|^4 d\tau$ is large, which implies that $E(W)$ is large.  We next turn to the case that $E(W)$ is very large.

Before moving on,  let us briefly describe this argument from the point of view of the tensor $S_{M_\Phi, r}$.   Recall that $S_{M_\Phi, r}$ is a multilinear map from $\CC^T \times ... \times \CC^T$ to $\CC$.   When we take the Fourier transform of $W$,  we are choosing a new orthonormal basis of $\CC^T$ and rewriting the tensor using that basis.   Equations (\ref{momsingr2}) or (\ref{momsingrint}) describe the tensor $S_{M_\Phi, r}$.   For any set of frequencies $\Phi$,  the tensor $S_{M_\Phi, r}$ in this basis is somewhat sparse,  because the difference set $D^r \Phi$ is supported on the hyperplane $\sum_{i = 1}^r \tau_i = 0$.   We also want to decompose $S_{M_\Phi,r} = N^r I + S^\Delta_{M_\Phi, r}$.   In the new basis,  both $N^r I$ and $S^\Delta_{M_\Phi, r}$ are supported on this hyperplane.   But for the special choice of $\Phi = \Phi_{Dir} =  \{ \log n \}_{n=N+1}^{2N}$,  the tensor $S^\Delta_{M_\Phi, r}$ is even sparser,  essentially supported on only a small subset of the hyperplane.




\subsection{The case with maximal additive structure} \label{subsecaddstruc}



When $W$ has a lot of additive structure, we get strong estimates using Heath-Brown's work \cite{HB} from the late 70s. 
This work can be seen as a refinement of the Hardy-Littlewood majorant principle, and we begin by reviewing that work.

The majorant principle is related to the following question.
Suppose that $D(t)$ is a trigonometric polynomial

\begin{equation} \label{Dtrigpoly} D(t) = \sum_{\xi \in \Phi} b_{\xi} e^{i t \xi}.
\end{equation}

\noindent Here $\Phi$ is a finite set and $b_\xi \in \CC$.  Suppose that for each $\xi$,  $B_\xi \ge |b_\xi|$.  (In particular, $B_\xi$ is real and non-negative.)   We want to compare $D(t)$ to the ``majorizing'' trigonometric polynomial

\newcommand{\Dmaj}{D_{\textrm{major}}}

\begin{equation} \label{Dtrigpoly} \Dmaj(t) = \sum_{\xi \in \Phi} B_{\xi} e^{i t \xi}.
\end{equation}

\noindent The general question raised by Hardy and Littlewood is,  ``how does the size of $D$ compare with the size of $\Dmaj$?''  This question is motivated by the following proposition.

\begin{prop} \label{propHLeven} Suppose that $\Phi \subset 2 \pi \ZZ$ so that $D(t)$ is a Fourier series on $\RR / \ZZ$,  and suppose $D$ and $\Dmaj$ are as above.   If $s$ is a natural number,  then

$$\int_0^1 |D(t)|^{2s} dt \le \int_0^1 | \Dmaj (t) |^{2s} dt. $$

\end{prop}

\begin{proof} The left hand side is

$$ \int_0^1 | D^s |^2 = \sum_{n \in \ZZ} \left| \widehat{D^s} \right|^2 = \sum_{n \in \ZZ} \left| \hat D^{*s}(n) \right|^2. $$

Similarly, the right-hand side is $ \sum_{n \in \ZZ} \left| \widehat{ \Dmaj}^{*s}(n) \right|^2. $

Here we write $\hat D * \hat D( n )$ for the convolution in $\ZZ$ and we write $\hat D^{*s}$ for $\hat D$ convolved with itself $s$ times.    

By our majorant hypothesis,  we have $|\hat D(n)| \le \widehat{ \Dmaj}(n)$ for all $n$,  and so $| \hat D^{*s}(n) | \le \widehat {\Dmaj}^{*s} (n)$.  
\end{proof}

This Proposition involves even integer exponents,  and it does not extend to other exponents.   See \cite{GR} for a good overview of what we know about other exponents.

Heath-Brown realized that,  instead of integrating over the circle,  we can get similar estimates if we sum over sets with a special structure: the structure of difference sets. 

\begin{prop} \label{propHBHL}  Suppose that $D(t)$ and $\Dmaj(t)$ are as above.   If $\mathcal{T} \subset \RR$ is a finite set and $s$ is a natural number, then

$$\sum_{t_1, t_2 \in \mathcal{T}} |D(t_1 - t_2)|^{2s} \le \sum_{t_1, t_2 \in \mathcal{T}} |\Dmaj(t_1 - t_2)|^{2s} .  $$

\end{prop} 

\begin{proof}  We do the proof first for $s=1$.   Plugging in the definition of $D$ and expanding,  we get

$$\sum_{t_1, t_2 \in \mathcal{T}} |D(t_1 - t_2)|^{2} = \sum_{t_1, t_2 \in \mathcal{T}} \left| \sum_{\xi \in \Phi} b_\xi e^{i (t_1 - t_2) \xi} \right|^2 =  \sum_{t_1, t_2 \in \mathcal{T}}  \sum_{\xi_1, \xi_2 \in \Phi} b_{\xi_1} \bar b_{\xi_2}  e^{i (t_1 - t_2) (\xi_1 - \xi_2)}. $$

Next we switch the order of summation and rewrite the expression as

$$ =  \sum_{\xi_1, \xi_2 \in \Phi} b_{\xi_1} \bar b_{\xi_2}    \sum_{t_1, t_2 \in \mathcal{T}}  e^{i (t_1 - t_2) (\xi_1 - \xi_2)} =   \sum_{\xi_1, \xi_2 \in \Phi} b_{\xi_1} \bar b_{\xi_2}  \left| \sum_{t \in \mathcal{T}} e^{i t (\xi_1 - \xi_2)} \right|^2.$$

Now since $|b_{\xi}| \le B_\xi$,  we can bound this by

$$ \le   \sum_{\xi_1, \xi_2 \in \Phi} B_{\xi_1} \bar B_{\xi_2}  \left| \sum_{t \in \mathcal{T}} e^{i t (\xi_1 - \xi_2)} \right|^2.$$

Now unwinding we see that this is $ \sum_{t_1, t_2 \in \mathcal{T}} |\Dmaj(t_1 - t_2)|^{2}. $

This finishes the proof for $s=1$.   If $s > 1$, we apply the argument above to $D^s$.

\end{proof}

Heath-Brown used this principle to study Dirichlet polynomials of the form $D(t) = \sum_{n=N}^{2N} b_n e^{i t \log n}$ with $|b_n| \le 1$.   
He gave very good estimates for sums of the form $ \sum_{t_1, t_2 \in \mathcal{T}} |D(t_1 - t_2)|^2$.   In this case,  $\Dmaj(t) = \sum_{n=N}^{2N} e^{i t \log n}$.   It is conjectured that $|\Dmaj(t)| \le T^{o(1)} N^{1/2}$ for $1 \le |t| \le T$,  and the triangle inequality gives $|\Dmaj(t)| \lesssim N$ if $|t| \le 1$.   Combining this conjecture and Proposition \ref{propHBHL} would give

\begin{conj} \label{conjDirsumdiffset} Suppose that $D(t) = \sum_{n=N}^{2N} b_n e^{i t \log n}$ with $|b_n| \le 1$.   Suppose that $\mathcal{T}$ is a 1-separated set in $[0,T]$.  Then

$$ \sum_{t_1, t_2 \in \mathcal{T}} |D(t_1 - t_2)|^{2} \le T^{o(1)} \left( |\mathcal{T}| N^2 + | \mathcal{T}|^2 N \right). $$

\end{conj}

\noindent This conjecture would be sharp if it is true.   The conjecture that $|\Dmaj(t)| \le T^{o(1)} N^{1/2}$ for $1 \le |t| \le T$ is a deep open problem which also appeared in the Halasz-Montgomery argument.   For a fixed $N$ and $T = N^{O(1)}$,  it can be checked by brute force in time $\Poly(N)$.   Assuming the bound on $|D_{major}|$  is true,  this would give a $\Poly(N)$ time proof for Conjecture \ref{conjDirsumdiffset}.   This is in contrast to the main conjectures in the large value problem,  where there is no plausible approach in the literature for a $\Poly(N)$ length proof.

Heath-Brown studied the term $ \sum_{t_1, t_2 \in \mathcal{T}} |\Dmaj(t_1 - t_2)|^{2}$ in a clever way,  involving induction on $N$, and he was able to get quite close to Conjecture \ref{conjDirsumdiffset}.   In particular, he was able to prove Conjecture \ref{conjDirsumdiffset} in the regime $N \le T \le N^{3/2}$.   






In the context of our main problem, we would like to estimate sums of the form $ \sum_{t \in W} |D(t)|^2$ where $W \subset \ZZ \cap [0, T]$.   When $W$ has a lot of additive structure, then Heath-Brown's theorem is helpful.  
To see the idea, suppose first that $W$ is an arithmetic progression of the form $W = \{ j \alpha \}_{j=-J}^{J}$.    In this case, we choose $\mathcal{T} = W$ so that $\mathcal{T} - \mathcal{T}$ covers each element of $W$ with multiplicity $\sim J \sim |W|$.  Then we get $ \sum_{t \in W} |D(t)|^2 \lesssim |W|^{-1} \sum_{t_1, t_2 \in\mathcal{T}} |D(t_1 - t_2)|^2$ and we can apply Heath-Brown's theorem.  More generally, if $E(W) \approx |W|^3$, then we can choose $\mathcal{T}$ with $|\mathcal{T}| \sim |W|$ so that $\mathcal{T} - \mathcal{T}$ covers most elements of $W$ about $|W|$ times, and then we proceed in the same way.  

In the case that $E(W) \approx |W|^3$, this actually leads to the sharp bound $|W| \lessapprox N^{2 - 2 \sigma}$, matching Conjecture \ref{conjMontgomery} (cf.  Lemma 1.7 in \cite{GM}).

( This method applies to a wide range of trigonometric polynomials,  but for any given choice of frequencies $\Phi = \{ \phi(n) \}_{n=N}^{2N}$,  we would need to estimate the majorant $\sum_{n} e^{ i t \phi(n)}$.   In particular,  the method applies to the almost counterexample.   Suppose we are able to prove good bounds for the majorant.   (Such bounds can be checked in time $\Poly(N)$ if they are true.)   With these good bounds,  the method shows that if $\max_n |b_n| \le 1$ and $W$ has a lot of additive energy,  then $\sum_n b_n e^{i t \phi(n)}$ cannot be large on $W$.   Comparing these bounds with the almost counterexample,  we see that the condition $\max_n |b_n| \le 1$ leads to much better control on an arithmetic progression $W$ than does the condition $\sum_n |b_n|^2$.)
 
 The methods in \cite{GM} gave some progress on the large value problem, but the results are still far from Conjecture \ref{conjMontgomery} or Conjecture \ref{conjMontgomerylq}.  In the next sections, we discuss some barriers to proving these conjectures.

\section{A barrier from computational complexity theory} \label{secbarrierlowdeg}

Looking for polynomial time algorithms to certify that a random matrix obeys certain bounds is part of the field of average case computational complexity.  Researchers in computer science have studied many problems of this flavor.  They have developed general approaches to these problems together with a conjectural picture of which tasks can be done in polynomial time and which cannot.  As explained in \cite{DHPT}, this theory makes some predictions about how much a polynomial time algorithm can prove about the large value problem for a random matrix.  In this section, we give an introduction to this area (aimed at a broad audience including harmonic analysts and number theorists).  Then we reflect on what these ideas might suggest about proving estimates for large values of Dirichlet polynomials.

Many problems in this field are phrased in terms of distinguishing a random object from an object with a planted structure.  One of the earliest and best studied examples is the planted clique problem.  In this problem, you are given a graph on $n$ vertices which was constructed in one of two ways.

\begin{itemize}

\item A random graph, where each edge is included with probability 1/2, independently.

\item A graph with a planted clique.  A size $k$ is chosen ahead of time.  First, randomly select $k$ of the $n$ vertices and include all the edges between them.  Then include each other edge with probability 1/2, independently.

\end{itemize}

It is straightforward to check that with high probability, the largest clique in the random graph has size $O(\log n)$.  The problem is interesting when $k$ is larger than this, for instance $k = n^\alpha$ for some $\alpha \in (0,1)$.  The graph with a planted clique clearly has a clique of size $k$, much larger than the clique in a typical random graph.

The planted clique problem has several variants.   In one variant,  we fix $n$ and $k = n^\alpha$ ahead of time.   Then we are given a graph.   With probability 1/2, the graph was a random graph, and with probabiilty 1/2 the graph has a planted clique of size $k = n^\alpha$.   We have to decide whether the graph was random or has a planted clique, and we need to be correct with high probability.   In another variant,  we are given a graph with a planted clique and we have to find the clique.  In a third variant,  we are given a random graph, and we have to produce a certificate (a proof) that the graph has no large clique.  

When $k$ is bigger than $n^{1/2}$, then all variants of the problem can be solved efficiently.  But when $k$ is $n^{1/2 - \epsilon}$, then there is no known polynomial time algorithm to solve any version of the planted clique problem.  The problem has been studied intensively since the 1990s.  

Some researchers suspect that when $k \le n^{1/2 - \epsilon}$, there is no polynomial time algorithm to solve (any version of) the planted clique problem.  Proving that there is no polynomial time algorithm for such a problem is far out of reach, and would represent a major milestone in computer science.   But computer scientists have shown that the problem cannot be solved in polynomial time by certain types of algorithms.   One of the most interesting results here says that when $k \le n^{1/2 - \epsilon}$,  the planted clique problem cannot be solved by the sum of squares method.  
The sum of squares method is a general strategy for proving bounds on constrained optimization problems, including all the types of problems we have considered here.  There is a very readable introduction to the sum of squares method in \cite{BS}.  

We will not try to describe the sum of squares method in detail here, but we just make a few comments.  The method applies to constrained optimization problems of the following kind.  We are allowed to choose a point $x \in \RR^N$ subject to some constraints of the form $p_j(x) \le 0$, where $p_j$ are polynomials.  Subject to these constraints, we would like to maximize $P(x)$, another polynomial.  Usually the polynomials $p_j$ and $P$ are of low degree.  The method is based on a theorem from semi-algebraic geometry called the positivestellensatz which describes when a subset of $\RR^N$ defined by polynomial inequalities is non-empty.  The positivestellensatz is a generalization of the Hilbert nullstellensatz, which describes when a subset of $\CC^N$ defined by polynomial equations is non-empty.  The positivestellensatz involves sums of squares, which is where the method gets its name.

The sum of squares method depends on a degree $D$.  For a constrained optimization problem of complexity $N$, the degree $D$ sum of squares method runs in time $N^{O(D)}$.  As $D$ increases, the method takes longer, but the bounds potentially improve.  The method outputs a correct upper bound for the optimization problem, together with a polynomial length proof of the bound.  The bound is not necessarily sharp, but it is always correct.  There are two key features of the sum of squares method.

\begin{itemize}

\item The method is very versatile.  It applies to a wide range of problems of the flavor considered in this survey.  And for a wide range of problems, it delivers the best bounds of any known polynomial time algorithm.

\item There are theorems about the limits of the sum of squares method.  For many problems, computer scientists know quite precisely what bounds can and cannot be proven by the sum of squares method.

\end{itemize}

For instance, the sum of squares method can solve the planted clique problem when $k \ge n^{1/2}$,  and Barak, Hopkins, Kelner, Kothari,  Moitra,  and Potechin \cite{BHKKMP} showed that, for any $\epsilon > 0$, the sum of squares method (with constant degree) cannot solve the planted clique problem if $k \le n^{1/2 - \epsilon}$.  

The planted clique problem was one of the first well studied problems in this direction, but the theory has expanded to address many problems, including some of the problems we have met in this survey.

As a second example, we consider the problem of norms of tensors.  Suppose that $T$ is a random $k$-tensor on $\RR^n$, where each of the $n^k$ entries is chosen independently from the standard Gaussian $N(0,1)$.  With high probabiilty, the norm of the tensor $T$ is $\approx n^{1/2}$.  The sum of squares method can prove a bound $\| T \| \lessapprox n^{k/4}$ (cf.  \cite{HSS}),  and the sum of squares method with degree $O(1)$ cannot prove a bound of the form $\| T \| \lessapprox n^{k/4 - \epsilon}$ (cf.  \cite{PR}).

(Incidentally, if $k$ is even, then the bound of $n^{k/4}$ can be proven by flattening $T$ to a $n^{k/2} \times n^{k/2}$ matrix $\tilde T$ and then estimating $\| \tilde T \|$, as we discussed in the last section. )

In this survey, we are interested in certifying large value estimates for random matrices.  Recall that the large value problem to the sparse singular value problem.   We do not yet know exactly what bounds the sum of squares method can certify in these problems.  The state of the art is in the paper \cite{DHPT}.   The upper bounds from the paper \cite{DHPT} that we discussed in Section \ref{secschatten} all follow from the sum of squares method.  Building on work of Mao and Wein (\cite{MW}), the paper also makes some conjectures on the limits of the sum of squares method, and even the limits of polynomial time methods in general.

Following \cite{DHPT}, let us focus on the question when it is possible to prove a bound that significantly improves on the operator norm method -- the basic bound that comes from estimating $\| M \|$.  Let us recall the problem, and then we can state the results and conjectures from \cite{DHPT}.  Fix an exponent $\alpha$ in the range $1 < \alpha < 2$.  Let $N \rightarrow \infty$ and $T = N^\alpha$.  Suppose $M$ is a  random $T \times N$ matrix, where each entry is independent, and the entries are chosen from the standard Gaussian $N(0,1)$.    Suppose that $v \in \RR^N$ is a vector with $\| v \|_\infty \le 1$ and $| (Mv)_j | \ge N^\sigma$ for all $j \in W$.  Our goal is to estimate $|W|$.  We can compute $\| M \| \lesssim T^{1/2}$ and this certifies the basic bound

\begin{equation} \label{basicran}
 |W| \lesssim T N^{1 - 2 \sigma} = N^{\alpha + 1 - 2 \sigma}. \end{equation}

When can we certify a stronger bound?  The main result of \cite{DHPT},  Theorem \ref{maindhpt},  implies that we can certify stronger bounds in a certain range of parameters:

\begin{theorem} (Theorem 1.5 in \cite{DHPT}) In the setup above,  with a random $T \times N$ matrix $M_{ran}$,  if $\sigma > 1 - \alpha /4$,  then there exists $c(\alpha, \sigma) > 0$ so that

$$ |W| \lesssim N^{\alpha +1 - 2 \sigma - c}. $$

Moroever, this bound can be certified by the sum of squares method with degree at most a constant depending on $\alpha, \sigma$.  

\end{theorem}

On the other hand,  in the regime $\sigma \le 1 - \alpha / 4$,  the paper \cite{DHPT} conjectures that the constant degree sum of squares method cannot certify bounds stronger than the basic bound (\ref{basicran}).   This conjecture is based on earlier work of Mao-Wein \cite{MW}.

\begin{conj} \label{conjsos} (\cite{DHPT}, 2nd to last paragraph on page 3) In the setup above,  with a random $T \times N$ matrix $M_{ran}$,  if $\sigma \le 1 - \alpha/4$ and $c>0$,   then the sum of squares method with degree O(1) cannot prove that

$$ |W| \lesssim N^{\alpha + 1 - 2 \sigma-c}. $$

\end{conj}

This conjecture is not yet proven,  but it is closely related to the theorems about planted clique and tensor norms mentioned above,  and it fits into a general program in the field.   The paper \cite{DHPT} also states an even stronger conjecture: if $\sigma < 1 - \alpha / 4$,  then there is no polynomial time algorithm that proves a large value estimate stronger than the basic bound (\ref{basicran}).   This latter conjecture is clearly out of reach of current methods,  but it is very intriguing, and it fits into a general picture about which tasks in this field can be solved in polynomial time.

Next we describe some of the ideas that motivate Conjecture \ref{conjsos}.   These ideas give a new point of view about why it is hard to improve on the basic bound (\ref{basicran}) when $\sigma \le 1 - \alpha/4$.   These ideas come from the theory of low degree testing.   We will explain this theory and state a rigorous theorem of Mao and Wein \cite{MW} which implies that proofs based on low degree tests cannot improve on (\ref{basicran}) when $\sigma < 1 - \alpha/4$.  

Our setup is similar to the planted clique problem.   We are given an $N^\alpha \times N$ matrix $M$ that was constructed in one of two ways.

\begin{itemize}

\item A random matrix,  where each entry is chosen from the standard Gaussian $N(0,1)$.

\item A matrix with a planted structure that makes the bound (\ref{basicran}) sharp.   The planted structure involves a sparse vector in the column space of $M$ -- we will describe this structure in detail below.

\end{itemize}

The problem is whether we can distinguish the two types of matrices (with high probability).   Here is a construction of matrices with planted structure that was studied by Mao-Wein in \cite{MW}.

\begin{enumerate}

\item The variables $N$, $ \alpha$, and $\sigma$ are given,  with $1 < \alpha < 2$ and $1/2 < \sigma < 1$.   Set $T = N^\alpha$.  

\item Set $S = N^{\alpha + 1 - 2 \sigma - \epsilon}$ for some tiny $\epsilon > 0$.

\item Randomly choose a subset $W \subset \{1, ..., T \}$ of cardinality $S$.

\item Randomly choose a sparse vector $w \in \RR^T$ supported in $W$.   If $j \in W$,  we choose $w_j$ from the Gaussian $N(0,  (T/S)^{1/2})$.  If $j \notin W$,  we set $w_j = 0$.

\item Build a $T \times N$ matrix $A$ where the first column of $A$ is $w$ and every other entry is chosen independently from the Gaussian $N(0,1)$.

\item Define $M_{plant} = A O$,  where $O$ is  random $N \times N$ orthogonal matrix.

\end{enumerate}

Because of the planted structure in this matrix,  the bound (\ref{basicran}) is essentially sharp.   We set $e_1 \in \RR^n$ to be the first component vector and we set $v = O^{-1}( \sqrt{N} e_1 )$.   Then with high probability,  each entry of $v$ has norm $\lessapprox 1$.   Note that $M_{plant} v = \sqrt{N} w$,  and so for most $j \in W$,  we have

$$ | (M_{plant} v)_j | \gtrapprox N^{1/2} (T/S)^{1/2} = N^\sigma. $$

This planted problem helps us reflect on improving the bound (\ref{basicran}),  because any proof that improves the bound (\ref{basicran}) for a random matrix $M_{ran}$ must somehow distinguish the matrix $M_{ran}$ from $M_{plant}$. 

Many of the proofs we have seen are based on bounded degree polynomials in the entries of the matrix.   For instance,  consider the approach from \cite{DHPT} based on the Schatten tensor.   The Schatten tensor $S_{M, r}$ has entries that are polynomials in the entries of $M$.    Following \cite{DHPT} (or \cite{GM}),  we may decompose this tensor into a few other tensors,  call them $T_M$.    The entries of $T_M$ will still be polynomials in entries of $M$.   We flatten each $T_M$ to get a matrix $\tilde T_M$.   Next we need to estimate the norm $\| \tilde T_M \|$.   One natural way to do this is to look at $\tr  \left( ( \tilde T_M \tilde T_M^*)^\ell \right)$ for some natural number $\ell$.   Expanding out,  this trace is a polynomial in the entries of $M$.   Let us call the final polynomial $p$.    Whenever the method in \cite{DHPT} improves on the basic bound (\ref{basicran}),  this polynomial $p$ distinguishes random matrices $M_{ran}$ from planted matrices $M_{plant}$.    The polynomial $p$ is always non-negative,  $p(M_{ran})$ is small with high probability, and $p(M_{plant})$ is large with high probability.    The degree of $p$ is a large constant (independent of $N$). 

Low degree testing is a general method for trying to distinguish between a random object and an object with planted structure.   Roughly speaking,  low degree testing
asks whether there is any low degree polynomial that distinguishes $M_{ran}$ from $M_{plant}$ - we will give the precise statement below.   In the paper \cite{MW},  Mao and Wein completely analyzed low degree testing for the random and planted matrix models described above.  We need a little notation to say precisely what they proved.   Suppose $\mu_{ran}$ is the probability distribution for random matrices $M_{ran}$ above and $\mu_{plant}$ is the probablity distribution for planted matrices $M_{plant}$ above.    We say that a polynomial $p$ is normalized if

\begin{equation} \label{lowdegnorm}
\int p(M) d \mu_{ran}(M) = 0 \textrm{ and } \int |p(M)|^2 d \mu_{ran}(M) = 1.
\end{equation}

Low degree testing asking if there is any normalized polynomial $p$ of degree at most $D$ so that

\begin{equation} \label{lowdegtest}
 \int p(M) d \mu_{plant}(M) > 2.
\end{equation}

Mao and Wein analyzed this question precisely, proving

\begin{theorem} (\cite{MW}, Theorem 4.5) Suppose that $N, \alpha, \sigma$ are given with with $1 < \alpha < 2$ and $1/2 < \sigma < 1$.    Let $M_{ran}$ and $M_{plant}$ be random matrices constructed as above.

If $\sigma > 1 - \alpha/4$,  then there is a normalized polynomial $p$ obeying (\ref{lowdegtest}) with degree at most $D(\alpha, \sigma)$.

If $\sigma < 1 - \alpha/4$,  then there is no normalized polynomial $p$ obeying (\ref{lowdegtest}) with degree at most $(\log N)^{O(1)}$.  

\end{theorem}

For the reader who is not familiar with this subject, the precise definition of low degree testing may look a little strange.   One might consider formulating low degree testing in other ways.  For instance,  one might ask if there is a polynomial $p$ of degree at most $D$ and a constant $c$ so that $\Prob_{\mu_{ran}} [ p(M) > c]$ is very different from $\Prob_{\mu_{plant}} [ p(M) > c]$.   This latter question is still open and would be interesting to understand.    There are several motivations for the formulation of low degree testing given here.   For this survey,  we mention two of them.

\begin{itemize}

\item This formulation of low degree testing is approachable.   It has been completely analyzed for a broad range of problems about planted structure,  and the proofs are fairly clean.

\item For a broad range of problems,  the best low degree tests correspond to the best known polynomial time algorithms.  

\end{itemize}

In \cite{Ho},  Hopkins formulated low degree testing and he observed that for a broad class of problems,  low degree testing tracks the best polynomial time algorithms.   For instance,  this is true for the planted clique problem and the tensor norm problem.   For planted clique,  if $k \ge n^{1/2}$, then low degree testing succeeds,  and there is a polynomial time algorithm to decide whether a graph has a planted clique.  If $k < n^{1/2 - \epsilon}$, then low degree testing fails and there is no known polynomial time algorithm to do so.   The situation is similar for many different problems about planted structures. 

In \cite{Ho},  Hopkins formulated a conjecture saying that for a precise class of problems,  low degree testing tracks the best polynomial time algorithms.   This conjecture was slightly revised by Holmgren and Wein in \cite{HW}.   

Low degree testing grew out of work on the sum of squares method such as \cite{BHKKMP}.  For many problems in the broad class,  in the regime where low degree testing fails,  researchers have proved that the sum of squares method also fails.   For instance,  this holds for the planted clique problem and the tensor norm problem.   This connection between low degree testing and sum of squares motivates Conjecture \ref{conjsos}.

To summarize,  in the setup above where $M_{ran}$ is a random $T \times N$ matrix,  in the regime $\sigma < 1 - \alpha / 4$,  there may well be no polynomial time method to certify large value estimates better than the basic bound (\ref{basicran}).   There may even be no polynomial length proofs of large value estimates better than (\ref{basicran}).   If there are polynomial length proofs or polynomial time methods to prove such bounds,  it would be a striking development in average case computational complexity and it might well lead to better algorithms for many problems in the subject.   Such a proof could not be based on low degree testing, and
if Conjecture \ref{conjsos} is correct,  such a proof could not be a sum of squares proof.    The sum of squares method is flexible enough to prove every bound for the large value problem for a random matrix $M$ that we have mentioned,  and so the proof would have to look very different from the proofs we have considered.

\subsection{Implications for Dirichlet polynomials?}

In the last subsection,  we have described a significant barrier to proving large value estimates for random matrices in the regime $\sigma < 1 - \alpha/4$.   (Recall that we are considering random $T \times N$ matrices where $T = N^\alpha$. )  In this regime,  there may well be no polynomial length proofs that go beyond the basic bound (\ref{basicran}).   And I think it is very likely that there is no sum of squares proof of degree $O(1)$ that goes beyond this basic bound.

Does this barrier have any implications for proving large value estimates for Dirichlet polynomials?  Currently,  in the regime where $\sigma < 1 - \alpha/4$,  we have no large value estimates beyond the basic bound (\ref{basicorthintro}),  which matches (\ref{basicran}).   For instance,  if $\alpha = 6/5$,  we have no new large value estimates when $\sigma < 1 - \alpha/4 = 7/10$.

If an explicit matrix $M$ seems to be pseudo-random,  then it might be reasonable to conjecture that it behaves similarly to a random matrix in terms of what large value estimates can be proven by the degree $D$ sum of squares method.   For instance,  if $M_\pi$ is a $T \times N$ matrix with $\pm 1$ entries corresponding to the digits of $\pi$ in base 2,  then I think it is very reasonable to conjecture that the degree $D$ sum of squares method proves the same large value estimates for $M_\pi$ that it does for $M_{ran}$.   Another candidate matrix might be defined by setting $M_{t,n} = e^{2 \pi i P(t,n)}$ where $P(t,n)$ is a messy polynomial such as

$$ P(t,n) = \sqrt{2} t^5 n^3 + \sqrt{3} t^2 n^6 + \sqrt{7} t^4 n^4 + \sqrt{11} t^2 n^6. $$

\noindent Let us write $M_P$ for this $T \times N$ matrix.   I think it is a reasonable conjecture that the degree $D$ sum of squares method proves the same large value estimates for $M_P$ that it does for a random matrix $M_{ran}$.  

This discussion raises a philosophical question: to what extent does $M_{Dir}$ behave like a random matrix $M_{ran}$?  We have seen several ways in which $M_{Dir}$ does not behave like a random matrix.   For instance,  the approximate geometric series example from Section 
\ref{subsecappgeo} is a non-random behavior,  and it shows that $M_{Dir}$ does not obey the same large value estimates for $\ell^2$ normalized input vectors as a random matrix (in other words, $M_{Dir}$ does not obey the conclusion of Proposition \ref{proprandlv}).   We have also seen that in an appropriate basis,  $S_{M_{Dir}, r}^\Delta$ is surprisingly sparse.   And in the case that $W$ has additive structure,  then $M_{Dir}$ interacts with it in a special way.

Some matrices from Fourier analysis have special structures that have been used to prove sharp bounds for $p$-to-$q$ norms or sharp large value estimates.   In some cases,  these estimates are stronger than estimates that we can certify for a random matrix in time $\Poly(N)$.   In Section \ref{secharmanal},  we will survey some of these structures. 

But the special structures for $M_{Dir}$ do not seem as useful or as important as the special structures that appear in these other matrices in Fourier analysis.    The matrix $M_{Dir}$ does have some special structures,  but the known ways of exploiting these structures do not lead to new estimates in the regime $\sigma < 1 - \alpha/4$.  

At this moment,  it may be interesting to return to the question whether we should believe Conjecture \ref{conjMontgomery} at all.   I personally do believe Conjecture \ref{conjMontgomery},  and my intuition is heavily guided by a pseudorandom argument,  which goes like this:

\begin{enumerate}

\item By Proposition \ref{proprandlv},  Conjecture \ref{conjMontgomery} holds for random matrices.   

\item The matrix $M_{Dir}$ has some special features that distinguish it from a random matrix.   For example,  the approximate geometric series example shows that the large value problem with $\ell^2$ normalized input behaves very differently for $M_{Dir}$ than for a random matrix.   However,  the known special features don't suggest any simple counterexample to Conjecture \ref{conjMontgomery}.

\end{enumerate}

\noindent Based on 1) and 2),  I suspect that Conjecture \ref{conjMontgomery} is true.   This type of pseudorandom argument is often used to make conjectures in analytic number theory.

Does this same type of pseudorandom argument extend to questions about whether there is a certain type of proof that $M_{Dir}$ obyes a certain bound?  I am not so sure,  but I think it is worth reflecting on.   And I am willing to make a conjecture that Conjecture \ref{conjMontgomery} cannot be proven by a sum of squares proof of degree $O(1)$.   In fact,  I am willing to conjecture that for each $\alpha \in (1,2)$,  there is some critical $\sigma > 1/2$ so that a sum of squares proof of degree $O(1)$ cannot improve on the basic bound (\ref{basicorthintro}).

\begin{conj} \label{conjmontsos} Consider the large value problem for $M_{Dir}$.   Suppose that $1 < \alpha < 2$ and set $T = N^\alpha$.   There exists $\sigma_{crit} (\alpha) > 1/2$ so that if $\sigma \le \sigma_{crit}(\alpha)$,  then for any degree $D$ and any $\epsilon > 0$ and any constant $C$,  for all $N$ sufficiently large,  there is no degree $D$ sum of squares proof that 

$$|W| \le C N^{\alpha + 1 - 2 \sigma - \epsilon}. $$

\end{conj}

One might even speculate whether $\sigma_{crit}(\alpha) = 1 - \alpha/4$,  which would match the numerology in Conjecture \ref{conjsos} for random matrices.

Many proofs in the literature on the large value problem can be rewritten as sum of squares proofs.   The sum of squares method can definitely do the basic operator norm method,  the power trick, and the $M M^*$method -- see Fact 1.3 and Fact 1.4 in \cite{DHPT}.   I'm not sure whether \cite{HB} or \cite{GM} can be rewritten as sum of squares proofs of bounded degree.   It would be interesting to check this carefully.  I would guess that all these proofs can be rewritten as sum of squares proofs of bounded degree.

In this problem, and many old problems in analytic number theory,  people who have worked on them have a feeling that ``we can't solve this problem with the methods we've been using,  and a proof of the full conjecture should look different."  It would be interesting to make (and prove) a precise statement that proofs of a certain kind cannot prove a certain result.   In analytic number theory and harmonic analysis,  we have rarely succeeded in doing this.  (One notable example is Selberg's work on the parity problem and the twin primes conjecture - see \cite{T2} for an introduction.   Another example is Katz's analysis of a certain approach to the Kakeya problem in \cite{K}.)  I spent some time trying to precisely define some limited system of proofs that would include most of the known proof methods for the large value problem,  but I wasn't happy with any definition I came up with -- they all seemed very messy and somewhat arbitrary.   The sum of squares method of degree $D$ is a precisely defined proof system that has been succesfully studied for many different problems in computer science,  and it may also be natural to consider in this context.   In particular,  I think it would be interesting to know what proofs in the harmonic analysis literature can be converted to sum of squares proofs and which cannot. 

To summarize,  it is not completely clear how much the low degree testing barrier for random matrices has to say about $M_{Dir}$.   The ideas from computer science do raise the question whether Conjecture \ref{conjmontsos} holds,  and perhaps they give some evidence that it does hold.   Conjecture \ref{conjmontsos} implies that there is no constant degree sum of squares proof of Conjecture \ref{conjMontgomery} or Conjecture \ref{conjMontgomerylq}.   A proof of either conjecture would have to be complex in a certain precise sense.   Perhaps that might even give a little clue of what we have to look for to make more progress.

\section{A barrier related to the Kakeya problem} \label{secbourkak}

In \cite{Bour},  Bourgain found a connection between the large value problem for Dirichlet polynomials and the Kakeya problem in geometric measure theory / harmonic analysis.   We will explain this connection here, and we will recall the Kakeya problem when it comes up.   This barrier is related to some special structure of the matrix $M_{Dir}$,  and it wouldn't come up for a random matrix.   The special structure is connected with the almost geometric progressions in Section \ref{subsecappgeo}. 

As we discussed in Section \ref{subsecappgeo},  short Dirichlet polynomials can approximate geometric series.   If $I \subset \{ N+1, ..., 2N \}$,  is an interval,  we write

$$D_{I, 0}(t) := \sum_{n \in I} e^{ i t \log n}. $$

\noindent If the interval $I$ is short,  then $D_{I,0}$ is close to a geometric series,  and as a result $D_{I,0}$ is large at many points.    The set where $D_{I, 0}$ is large is called a fat arithmetic progression,  which means the $r$-neighborhood of an arithmetic progression $A$,  for some $r$ and $A$.  
For each $I$,  there is a fat arithmetic progression $I^*$ so that

$$ | D_{I,0}(t) | \sim |I| \textrm{ for all } t \in I^*. $$

As $I$ gets longer,  $D_{I,0}$ behaves less like a geometric series,  and as a result the range of $t$ where the geometric series approximation holds gets smaller.   If we want to study Dirichlet polynomials on $[0,T]$,  then there is a critical length of $I$ so that $D_{I,0}$ is well modeled by a geometric series on all of $[0,T]$.  The critical length works out to be $|I| = N T^{-1/2}$.   For this length,  we have that for all $|t| \le T$,

$$ | D_{I,0}(t) | \sim |I| 1_{I^*}(t). $$


We mentioned that these almost geometric series are counterexamples to some $\ell^2$ based large value conjectures.   But they are not counterexamples to the main conjectures because in this example,  the coefficients $b_n$ are non-zero only on a small interval $I \subset \{ N, ..., 2N\}$,  whereas the conjectures are about Dirichlet polynomials with $|b_n| \sim 1$ for all $n \in \{N, ..., 2N \}$.   We can build interesting examples of this type by combining $D_I$ for many different $I$.    To give ourselves more flexibility,  we will also use translations of the $D_{I,0}$.   

Next we recall that a translation of a Dirichlet polynomial is again a Dirichlet polynomial.   Suppose that $D(t) = \sum_{n \in I} b_n e^{i t \log n}$.    Define $D_{t_0} (t) = D(t - t_0)$.   Notice that $D_{t_0}$ is itself a Dirichlet polynomial:

$$ D_{t_0} (t)  = \sum_{n \in I} b_n e^{- i t_0 \log n} e^{i t \log n}. $$

\noindent This is a Dirichlet polynomial with coefficents $\tilde b_n = b_n e^{-i t_0 \log n}$.    

Now we can consider translations of $D_{I,0}$.   We let $D_{I, t_I}(t) = D_{I,0}(t-t_I)$.   We write $I^* + t_I$ to denote the translation of $I^*$ by $t_I$, and we note that $D_{I,t_I}$.  In particular, we can say that

$$ | D_{I, t_I}(t)| \gtrsim |I| 1_{I^* + t_I}. $$

\noindent Here $|I|$ is the cardinality of $I$,  and $1_X$ denotes the characteristic function of $X$.

Now we are ready to combine different $I$.   We decompose $\{ N, ..., 2N\}$ as a disjoint union of intervals $I$, and now we consider Dirichlet polynomials of the form

$$ D(t) = \sum_I \pm D_{I, t_I}(t). $$

\noindent In this setup,  for each $I$ we get to choose a translate $t_I$ and a sign $\pm 1$.   If we choose the signs randomly,  then with high probability we get

$$ |D(t)|^2 \approx \sum_I | D_{I, t_I}(t) |^2 = |I|^2 \sum_I 1_{I^* + t_I}, $$

and so for all $p > 2$, 

$$ \int_0^T |D(t)|^p dt \approx |I|^p \int_0^T \left( \sum_I  1_{I^* + t_I} \right)^{p/2}. $$

If the sets $I^* + t_I$ overlap heavily,  then we would get an interesting example for the large value problem.   Deciding whether it is possible to choose $t_I$ so that $I^* + t_I$ overlap heavily is a Kakeya type problem.

Let us pause here to recall the original Kakeya problem.
The original Kakeya problem concerns tubes pointing in different directions in $\RR^d$.   We divide $S^{d-1}$ into $1/R$-caps $\theta$.   For each $\theta$, we let $\theta^*$ denote a tube of radius 1 and length $R$,  pointing in the direction perpendicular to $\theta$, and centered at the origin.   For each $\theta$,  we consider a translate $\theta^* + x_\theta$.   This is a tube of radius 1 and length $R$ pointing perpendicular to $\theta$,  but now it may have any center.   The Kakeya problem is about deciding whether it is possible to choose $x_\theta$ so that the sets $\theta^* + x_\theta$ overlap heavily.   There are several different ways to make this precise.   One way is to try to maximize the $L^p$ norm $ \| \sum_\theta 1_{\theta^* + x_\theta} \|_{L^p}$ for different $p$.   A second way is to try to minimize the volume of the union $\cup_\theta (\theta^* + x_\theta)$.   

If we choose $x_\theta = 0$ for every $\theta$,  then we find that the tubes fill $B_R^d$ and the volume of the union is $R^d$.   Around 1920,  Besicovitch found a clever arrangement of tubes with volume $\sim (\log R)^{-1} R^d$.   But no one has even found an arrangement of volume $R^{d - \epsilon}$.   The Kakeya conjecture for volume asserts that there is no such arrangement.

\begin{conj} \label{kakvol} Fix a dimension $d \ge 2$ and let $\theta$ and $\theta^*$ be as above.  For every $\epsilon > 0$,  there is a constant $C_{\epsilon, d}$ so that for every choice of $x_\theta$,  

$$ \left| \bigcup_\theta (\theta^* + x_\theta) \right| \ge C_{\epsilon,d} R^{d - \epsilon}. $$

\end{conj}

There is a similar conjecture for $L^p$ norms.   If we choose all $x_\theta = 0$,  then it is straightforward to check that

$$ \int \left| \sum_\theta 1_{\theta^*} (x) \right|^p dx \approx R^{(d-1)p} + R^d. $$

\noindent The contribution $R^{(d-1)p}$ comes from the unit ball around 0,  which lies in all $R^{d-1}$ tubes.   The contribution $R^d$ comes because the integrand is 1 at most points of $B_R$.   The $L^p$ version of the Kakeya conjecture says that up to factors of the form $R^\epsilon$, the $L^p$ norm is always bounded by this example.   

\begin{conj} \label{kaklp} Fix a dimension $d \ge 2$ and let $\theta$ and $\theta^*$ be as above.   For every $\epsilon > 0$ there is a constant $C(\epsilon, d)$ so that for every $p \ge 1$ and every choice of $x_\theta$, 

$$ \int \left| \sum_\theta 1_{\theta^* + x_\theta} (x) \right|^p dx \le C(\epsilon, d) R^\epsilon \left( R^{(d-1)p} + R^d \right). $$

\end{conj}

These Kakeya conjectures about tubes have natural analogues in the setting of fat arithmetic progressions.   The $L^p$ version is the following.

\begin{conj} \label{bourkaklp} Fix $\alpha \in (1,2)$ and suppose that $T = N^\alpha$.    Divide $\{N+1, ..., 2N\}$ into intervals $I$ of length $\sim N T^{-1/2}$.   Let $I^*$ be as above.   There is a constant $C(\epsilon, \alpha)$ so that for every $N$ and every $p \ge 1$ and every choice of $t_I \in [0,T]$,  

\begin{equation} \label{conjBK}
 \int_0^T \left( \sum_I  1_{I^* + t_I} \right)^{p/2} \le C(\epsilon) N^\epsilon  \int_0^T \left( \sum_I  1_{I^*} \right)^{p/2}. 
\end{equation}

\end{conj}

An easy corollary of Conjecture \ref{bourkaklp} is that for any choice of $t_I$,

\begin{equation} \label{conjBKvol}  | \cup_I I^* + t_I | \ge C(\eps) N^{-\epsilon} T \approx | \cup_I I^*|. 
\end{equation}

By considering the special examples above,  Bourgain showed that Conjecture \ref{conjMontgomerylq} implies Conjecture \ref{bourkaklp} and hence
(\ref{conjBKvol}).   He also showed that these conjectures imply Conjecture \ref{kakvol},  the original Kakeya conjecture for tubes.  Because of Bourgain's work introducing this problem,  I propose to call this question the Bourgain-Kakeya problem and to call Conjecture \ref{bourkaklp} the Bourgain-Kakeya conjecture.  

When $d=2$,  both versions of the Kakeya conjecture (Conjecture \ref{kakvol} and Conjecture \ref{kaklp}) were proven by Davies around 1980.  
For $d \ge 3$,  both versions of the Kakeya problem were open for a long time.   Very recently in \cite{WZ},  Wang and Zahl announced a proof of Conjecture \ref{kakvol} for $d=3$.   Several related problems have also been solved,  including the finite field analogue (Dvir \cite{D}),  the p-adic analogue (Arsovski \cite{A} and Dhar-Dvir \cite{DD}),  and the Furstenberg set problem (Ren-Wang \cite{RW} and Orponen-Shmerkin \cite{OS}).  
It's not yet clear how much the methods in these works will say about the Bourgain-Kakeya problem.   

It is also worth mentioning that assuming the Bourgain-Kakeya conjecture as a black box does not currently lead to any improvements in the known bounds on the large value problem for Dirichlet polynomials.

\section{Examples from harmonic analysis} \label{secharmanal}

In harmonic analysis,  researchers have studied $p$-to-$q$ norm estimates and large value estimates for many different matrices $M$.   
There is a particularly close parallel between the large value problem for Dirichlet polynomials and the area of Fourier analysis called restriction theory.  Broadly speaking, restriction theory can be defined as the study of the matrices $M_\Phi$ for some set of frequencies $\Phi$.  In restriction theory, the set of frequencies $\Phi$ usually has a geometric structure, like a submanifold of $\RR^d$.

The classical methods from Section \ref{secclassmeth} are used widely in harmonic analysis and especially in restriction theory.  
The operator norm method is applied very broadly, where the operator norm of $M$ is often estimated using orthogonality, as in Plancherel's theorem.  The power method and the $M M^*$ method are also used broadly.  These general methods give estimates for many different problems.  Sometimes these estimates are sharp, but often they are not.  We will discuss some problems where the special structure of the matrix $M$ can be used to go beyond these methods.  

\subsection{Restriction theory}

Suppose that $S \subset \RR^d$ is a submanifold, such as the unit sphere or the paraboloid defined by $\omega_d = \sum_{j = 1}^{d-1} \omega_j^2$.  We let $d \mu_S$ be the volume measure on $S$.  We define the operator $E_S$ by

\begin{equation} \label{defEs}
E_S g (x) := \int e^{2 \pi i \omega x} g(\omega) d \mu_S(\omega). 
\end{equation}

Here $g$ is a function on $S$ and $E_S g$ is a function on $\RR^d$.  One motivation for this comes from linear PDE.  A solution to a given linear PDE can be written in the form $E_S$, where the submanifold $S$ encodes the PDE.  For example, if we think of $x_d$ as a time coordinate $t$, then a solution to the free Schrodinger equation $\partial_t u(x,t) = i \sum_{j=1}^{d-1} \partial_j^2 u$ can be written in the form $E_S$ where $S$ is the paraboloid $\omega_d = \sum_{j=1}^{d-1} \omega_j^2$.

In the late 60s, Stein raised the question of the $p$-to-$q$ norm of the linear operator $E_S$ for different surfaces $S$ beginning with the unit sphere.  The main question here is to determine all exponents $p,q$ so that

\begin{equation} \label{restr}  \| E_{S^{d-1}} g \|_{L^q(\RR^d)} \lesssim \| g \|_{L^p(S^{d-1})}. \end{equation}

\noindent Stein noticed that curved submanifolds behave differently from flat submanifolds and made an influential conjecture about the behavior of $\| E_S \|_{p \rightarrow q}$, called the restriction conjecture.

The linear operator $E_S$ is described in a continuous way, but we can discretize the situation to make a matrix $M_s$ which essentially captures the behavior of the linear operator $E_S$.  For example, suppose $S$ is the paraboloid.  Because of scaling symmetry, we can usually reduce questions to the case that $g$ is supported on $\{ \omega \in S: |\omega| \le 1 \}$.  Now fix a large parameter $R$.  In place of $E_S$, consider exponential sums of the form

\begin{equation} \label{trigpoly}
 F(x) = \sum_{n \in \ZZ^{d-1}, |n| \le R} b_n e^{2 \pi i ( \sum_{j=1}^{d-1} \frac{n_j}{R} x_j + \frac{ \sum_j n_j^2}{R^2} x_d ) } . 
\end{equation}

Finallly consider a matrix $M_S$ where the columns are indexed by $\{ n \in \ZZ^{d-1}, |n| \le R \}$ and the rows are indexed by $\{x \in \ZZ^d, |x| \le R \}$, with $M_{x,n} = e^{2 \pi i ( \sum_{j=1}^{d-1} \frac{n_j}{R} x_j + \frac{ \sum_j n_j^2}{R^2} x_d ) } . $ So $M_S$ is roughly a $R^d \times R^{d-1}$ matrix.  Note that for $x \in \ZZ^d, |x| \le R$, we have $F(x) = (M_S b)_x$.

Estimating $\| E_S \|_{p \rightarrow q}$ is essentially equivalent to estimating $\| M_S \|_{p \rightarrow q}$ uniformly in $R$ as $R \rightarrow \infty$.  One can perform a similar discretization for other surfaces $S$ such as the sphere.

The lecture notes \cite{T} are a good reference for the restriction problem and give details for everything in this subsection.  

The restriction problem was studied using the fundamental techniques we introduced in Section \ref{secclassmeth}.  In unpublished work in the late 60s, Stein studied $E_S$ using the $M M^*$ method, proving that $E_S$ obeys stronger estimates when $S$ is a sphere than when $S$ is a flat disk.  Around 1970, Fefferman \cite{F3} studied $E_S$ using the power method and proved the restriction conjecture in dimension $d=2$.  A few years later, Tomas and Stein added a wrinkle to the $M M^*$ method, and proved that $\| E_S \|_{2 \rightarrow q} \lesssim 1$ if and only if $q \ge \frac{2 (d+1)}{d-1}$.  

Up to this point, the methods for studying the restriction conjecture and the large value problem for Dirichlet polynomials are parallel.  In addition, the important examples are similar.  In the restriction conjecture, the relevant example is called the Knapp example.  It occurs when $g$ is a smooth supported on a small cap of the surface $S$.  In this case, the formula for $E_S g(x)$  is an approximate geometric series (or geometric integral), and the function $E_S g(x)$ is concentrated on a tube with long axis perpendicular to the tangent space of the small cap.  The approximate geometric series example in Section \ref{subsecappgeo} is analogous to the Knapp example.

Whenever a portion of $S$ is flat or almost flat, and $g$ is a smooth bump supported on that region, then $E_S g$ is an approximate geometric series.  It matters whether $S$ is curved or flat because that controls this type of example.  From now on, we suppose that $S$ is the unit sphere or the part of the paraboloid with $|\omega| \le 1$.  
  Because the Knapp example is supported on a small portion of $S$, one can hope for estimates of the form $\| E_S g \|_{L^q(\RR^d)} \lesssim \| g \|_{L^\infty(S)}$ with a wider range of $q$ than for estimates of the form $\| E_S g \|_{L^q(\RR^d)} \lesssim \| g \|_{L^2(S)}$.   When $S$ is the unit sphere or a compact piece of the paraboloid, then Tomas-Stein proved that $\| E_S g \|_{L^q(\RR^d)} \lesssim \| g \|_{L^2(S)}$ if and only if $q \ge \frac{2(d+1)}{d-1}$.  On the other hand, Stein conjectured that $\| E_S g \|_{L^q(\RR^d)} \lesssim \| g \|_{L^\infty(S)}$ if and only $q > \frac{2d}{d-1}$.  

Fefferman \cite{F3} proved the restriction conjecture in dimension $d=2$ using the power method, but the problem is still open in all dimensions $d \ge 3$.  Dimension $d=2$ is special because the critical exponent is $4$, which is an even integer.  In higher dimensions, the critical exponent is $\frac{2d}{d-1}$, which is not an even integer.

In other work, Fefferman showed that the restriction problem is related to the Kakeya problem.  In \cite{F1}, Fefferman used Besicovitch's construction of Kakeya sets with small volume to give a counterexample to a conjecture in Fourier analysis called the ball multiplier conjecture.  By similar reasoning, the restriction conjecture implies the Kakeya conjecture (Conjecture \ref{kaklp}).  This conjecture helped suggest Bourgain's work in the last section relating the large value problem for Dirichlet polynomials and the Kakeya conjecture.  

By 1975, people in restriction theory had used the operator norm method, the power method, and the $M M^*$ method.  In dimension $d=2$, these methods resolve the restriction problem, but in dimension $d \ge 3$ they do not.  It was open for about 15 years to find a bound going beyond these methods.  In the early 90s, Bourgain \cite{B3} proved estimates for (\ref{restr}) where $q = \infty$ and where $q$ is smaller than the Tomas-Stein range (i.e. $ q< \frac{2(d+1)}{d-1}$).  

Bourgain's proof is based on the idea of wave packets and on new bounds for the Kakeya problem.   Wave packets are the key structural tool that has been used to study $M_S$ over the last 30 years.   In Section \ref{subsecwavepack} below,  we will discuss this tool in more detail. 

Let us compare the restriction conjecture about the matrix $M_S$ to the bounds for random matrices.   The sharp bounds for $\| M_S \|_{2 \rightarrow q}$ (proven by Tomas-Stein) are not as strong as the bounds for $\| M_{ran} \|_{2 \rightarrow q}$  in Proposition \ref{proprandlv}.   The sharp examples for $\| M_S \|_{2 \rightarrow q}$ are the Knapp examples.   These Knapp examples constitute a special structure for $M_S$ that has no analogue for random matrices.   On the other hand,  the restriction conjecture for $\| M_S \|_{\infty \rightarrow q}$ matches the bounds for $\| M_{ran} \|_{\infty \rightarrow q}$ in Proposition \ref{proprandlv}. 

As we discussed is Section \ref{secbarrierlowdeg},  proofs based on low degree testing cannot prove optimal large value estimates for $M_{ran}$ in a certain regime.   In particular,  if $M_{ran}$ is a random $N^\alpha \times N$ matrix with $1 < \alpha < 2$,  then a proof based on low degree testing cannot improve on the simple bound $|W| \lessapprox N^{\alpha + 1 - 2 \sigma}$ coming from the operator norm method.   This implies that when $q > 6 - \frac{4}{\alpha}$,  then proofs based on low degree testing cannot prove bounds for $\| M_{ran} \|_{p \rightarrow q}$ that are sharp up to a factor $\lessapprox 1$.   For the matrix $M_S$ where $S$ is the unit sphere $S^{d-1} \subset \RR^d$,  the matrix $M_S$ has dimensions $R^d \times R^{d-1}$,  and so $\alpha = \frac{d}{d-1}$ and $6 - \frac{4}{\alpha} = 2 + \frac{4}{d}$.   While the restriction conjecture is still open,  state of the art results do give sharp bounds for $\| M_S \|_{\infty \rightarrow q}$ for $q < 2 + \frac{4}{d}$ -- see \cite{BG},  \cite{HR} and \cite{HZ}.

We can prove bounds for $\| M_S \|_{\infty \rightarrow q}$ that are out of reach to prove for a random matrix $M_{ran}$ because of the special structure of wave packets that we discuss more below.

\subsection{Large value estimates in restriction theory}

Knowing a bound on $\| M \|_{p \rightarrow q}$ gives bounds for the large value problem for $M$ with input normalized in $\ell^p$.  For a random matrix $M_{ran}$, the optimal bounds for the large value problem are stronger than the bounds coming from knowing $\| M_{ran} \|_{p \rightarrow q}$ for all $q$.  But in the restriction problem, the solution to the large value problem with input normalized in $\ell^2$ follows from Tomas-Stein.  And the conjectural solution to the large value problem with input normalized in $\ell^\infty$ would follow from the Stein restriction conjecture.  The sets $W$ that appear in the large value problem are small balls (or sometimes tubes).  We can make a more interesting large value problem by adding a spacing condition on the set $W$,  which says that $W$ cannot put too much mass into a small ball.  

In the 80s and 90s, Mattila and Wolff introduced and explored a large value question for $E_{S^1}$ with a stronger spacing condition.  This stronger spacing condition is related to Hausdorff dimension, and their goal was to study geometric properties of sets of a given Hausdorff dimension.  Here we rephrase their question a little bit to make a parallel with the large value problem.

Question: Suppose that $W \subset B(0, R) \subset \RR^2$ and that $| E_{S^1} g (x) | > \lambda$ for $x \in W$.  We make the normalization that $\| g \|_{L^2(S^1)} = R$.  Finally, suppose that $W$ obeys the following Frostman-type spacing condition: for any $B(x,r) \subset B(0, R)$, 

$$  | W \cap B(x,r) | \lesssim |W|^{\frac{\log r}{\log R}}. $$

\noindent What is the maximum possible size of $W$ in terms of $R$ and $\lambda$?

In \cite{Ma}, Mattila used the $M M^*$ method to show that if $\lambda \ge R^{3/4 + \eps}$, then $|W| \lesssim R^2 \lambda^{-2}$, which is tight.  His method did not give any bound for $\lambda \le R^{3/4}$.  Later Wolff \cite{W} proved a sharp bound for $\lambda \le R^{3/4}$: $W| \lessapprox R^4 \lambda^{-4}$.   Notice that there is a jump when $\lambda$ crosses $R^{3/4}$.  If $\lambda$ is just above $R^{3/4}$, then the $M M^*$ method gives $|W| \lesssim R^{1/2}$.  But if $\lambda$ is just below $R^{3/4}$, then the best estimate is $|W| \lessapprox R$.  

Wolff's bound can be proven either using the power method or using wave packets.

Wolff's bound is tight for all $\lambda$from $R^{1/2}$ up to  a small constant times $R^{3/4}$.   In particular, there is an example where $\lambda \sim R^{3/4}$ and $|W| \sim R$.   In this case, $g$ is supported on an arc of the circle of length $\sim R^{-1/2}$, and so $E_\mu g$ has size $\sim R^{3/4}$ on a rectangle of dimensions $R^{1/2} \times R$.  Then we take $W$ to be a set of $R$ points inside the rectangle, with $R^{1/2}$ points in each ball of radius $R^{1/2}$ in the rectangle.  (We cannot take any more points inside the rectangle because of the Frostman-type spacing condition.)   This example is based on the Knapp example and it is also similar in spirit to our almost counterexample.

In the last section of \cite{W}, Wolff raised a variant of the question above where the condition $\| g \|_{L^2(S^1)} = R$ is replaced by $\| g \|_{L^\infty(S^1)} = R$.   One may hope for a stronger upper bound on $|W|$ in this situation, since it rules out the Knapp example above.   Wolff explained how this problem is related to the Furstenberg set problem.   There hasn't been any progress on this question since Wolff's work,  although the recent breakthrough on the Furstenberg set problem (cf. \cite{RW}) probably leads to improved estimates here. 

There are a number of open problems related to this question.  For instance, the higher dimensional version of the questions in this section are largely open -- see \cite{DZ} for the best known estimates.

\subsection{The Strichartz estimate for the periodic Schrodinger equation}

A solution to the Schrodinger equation on the $(d-1)$-dimensional unit cube torus is a $\ZZ^d$-periodic function $u(x,t)$, where $x \in \RR^{d-1}$ and $t \in \RR$, which obeys the differential equation

$$ \partial_t u = - \triangle_x u := - (\sum_{j=1}^{d-1} \partial_{x_j}^2 u). $$

\noindent It is straightforward to check that such a solution $u$ is given by a Fourier series of the form

$$ u(x,t) = \sum_{n \in \ZZ^{d-1}} b_n e^{2 \pi i (n \cdot x + |n|^2 t)}. $$

We say that $u$ has frequencies up to $N$ if $b_n = 0$ for $|n| > N$.   Suppose that $u$ has frequencies up to $N$.   In the context of linear and non-linear PDE,  one is interested in estimates of the form

$$ \| u \|_{L^q([0,1]^{d})} \le C(q,N) \| b_n \|_{\ell^2}. $$

We can discretize this problem so that it becomes an estimate for $\| M_{PS} \|_{2 \rightarrow q}$ for an appropriate matrix $M_{PS}$.   To do this,  we note that if $\hat u(\xi)$ is supported on a rectangule of the form $|\xi_j| \le N$ for $j = 1, ..., d-1$ and $|\xi_d| \le N^2$.  Therefore,  the function $u$ is essentially constant on rectangles of dimensions $\frac{1}{N} \times ... \times \frac{1}{N} \times \frac{1}{N^2}$.   Therefore,  we define $M_{PS}$ to be a matrix where the columns are indexed by $n \in \ZZ^{d-1}, |n| \le N$ and the rows are indexed by $(x,t) \in [0,1]^d$ where $x \in (\frac{\ZZ}{N})^{d-1}$ and $t \in \frac{\ZZ}{N^2}$.    The entry $(M_{PS})_{(x,t), n}$ is given by $e^{2 \pi i (n \cdot x + |n|^2 t)}$ so that $(Mb)_{(x,t)} = u(x,t)$.  
The matrix $M_{PS}$ is a $N^{d+1} \times N^{d-1}$ matrix.   Bounding the quantity $C(q,N)$ above boils down to estimating $\| M_{PS} \|_{2 \rightarrow q}$.

In the early 90s, in\cite{B4},  Bourgain studied this problem.  He conjectured that

\begin{equation} \label{perstrich} C(\frac{2(d+1)}{d-1},  N) \le C_\eps N^\eps.  \end{equation}

This conjecture is equivalent to saying that $\| M_{PS} \|_{2 \rightarrow q}$ obeys the same bounds as $\| M_{ran} \|_{2 \rightarrow q}$ from Proposition \ref{proprandlv}.   We can also rephrase this estimate as a large value estimate.  
If we normalize so that $\sum_{|n| \le N} |b_n|^2 \lesssim N^{d-1}$,  and define $W_\lambda(u) := \{ (x,t): |u(x,t)| > \lambda$,  then this conjecture is equivalent to

\begin{equation} \label{perstrichlevelsets}  | W_\lambda(u) | \lessapprox N^{d+1} \lambda^{- \frac{2(d+1)}{d-1}}. 
\end{equation}

Note that the relevant range of $\lambda$ goes from $N^{\frac{d-1}{2}}$ to $N^{d-1}$.   The parameter $N$ in the context of Dirichlet polynomials corresponds to $N^{d-1}$ here.  

In this problem,  the critical exponent $q$ is $\frac{2(d+1)}{d-1}$.   This critical $q$ is an even integer if $d=2$ or 3,  but not for $d \ge 4$.   Using a version of the power method,  Bourgain proved in \cite{B4} that (\ref{perstrich}) holds when $d = 2$ or $3$.

Using a variant of the $M M^*$ method with an extra wrinkle,  Bourgain proved that (\ref{perstrichlevelsets}) holds for $\lambda >  N^{\frac{3}{4}( d-1)}$ (cf. Proposition 3.82 in \cite{B4}).   As we have seen in other scenarios,  the $M M^*$ method gives no information for $\lambda$ less than $N^{\frac{3}{4}(d-1)}$.  

As we have seen several times,  the operator norm method,  the power method,  and the $M M^*$ method give some bounds for a problem,  but the bounds do not appear to be sharp except in special cases ($d=2$ or 3).   

Much later,  in \cite{BD},  Bourgain and Demeter proved a major conjecture called the sharp $\ell^2$ decoupling conjecture.   This result implies in particular that (\ref{perstrichlevelsets}) holds for the whole range of $\lambda$ and for every $d$.   This upper bound for $|W_\lambda(f)|$ is sharp up to factors $\lessapprox 1$ for every $\lambda$  from $N^{\frac{d-1}{2}}$ to $N^{d-1}$.   In this case,  it turns out that there is no jump when $\lambda$ crosses $N^{\frac{3}{4}(d-1)}$. 

This theorem gives (\ref{perstrich}) and proves that $\| M_{PS} \|_{2 \rightarrow q}$ obeys the same bounds as $\| M_{ran} \|_{2 \rightarrow q}$ for all $q$. 

This theorem of Bourgain-Demeter was the first time we could prove such bounds for any specific $T \times N$ matrix $M$ with entries of norm 1 where $T$ is not an integer power of $N$.   Up to the present,  the only specific matrices for which we can prove such strong bounds are close cousins of $M_{PS}$ and the proofs are all based on decoupling.  The proof of decoupling is itself based on wave packet methods, and it takes advantage of wave packet structure at many different scales.  So the proof of these strong bounds for $\| M_{PS} \|_{2 \rightarrow q}$ is based on the special structure of $M_{PS}$ coming from wave packets.  We will discuss this a little more in Section \ref{subsecwavepack} below.

We mentioned in the introduction that there is no specific matrix $M$ for which we can prove the analogue of the sharp large value estimate which holds for random matrices with high probability (cf. Proposition \ref{proprandlv}).  For the matrix $M_{PS}$, this strong large value estimate does not hold.  There are a couple of different counterexamples, but one of them is a cousin of the Knapp example which is closely related to wave packets.

\subsection{Wave packets} \label{subsecwavepack}

The idea of wave packets has been central in restriction theory since early on, and especially since Bourgain's work \cite{B3}.  Wave packets are closely related to the Knapp example above.  So far, we have described how the Knapp example can be used to show that $\| E_S \|_{p \rightarrow q}$ is large.   It might be surprising that these types of examples can also be used to show that $ \| E_S \|_{p \rightarrow q}$ is small.  There are several good sources that explain this story in full, including Tao's notes \cite{T} and Demeter's book \cite{De}.  In this short section we will try to give a very non-technical overview.  We will also try, in our general context, to explain what is special about the matrix $M_S$ compared to a random matrix $M_{ran}$ or the Dirichlet matrix $M_{Dir}$.  

Suppose that we want to study $E_S g$ on a ball of some large radius $R$.  We subdivide $S$ into small caps $\theta$ of radius $\sim R^{-1/2}$.  Knapp noted that if $g$ is a smooth bump on $\theta$, then $E_S g |_{B_R}$ is concentrated in a tube of radius $R^{1/2}$ and length $R$, oriented perpendicular to the tangent space of $\theta$.  We call this tube $\theta^*$.  More generally, for any function $g_\theta$ supported on $\theta$, $E_S g_\theta |_{B_R}$ can be decomposed into pieces that are supported on translates of $\theta^*$ -- and each piece is called a wave packet.  Let us set up a little notation to describe this.

Suppose that $u = E_S g$.  For each cap $\theta$ in $S$, define $u_\theta = E_S (1_\theta g)$, so that $u = \sum_{\theta} u_\theta$.  Next let $\TT_\theta$ denote a tiling of $B_R$ by disjoint translates of $\theta^*$.  Then we can write

\begin{equation} \label{wptheta} u_\theta |_{B_R} = E_S g_\theta |_{B_R} = \sum_{T \in \TT_\theta} a_T W_T, \end{equation}

\noindent where $a_T \in \CC$ and $W_T$ is a function essentially supported on $T$, which is a translate of the Knapp example.  Each function $W_T$ is called a wave packet.  Since $u = \sum_\theta u_\theta$, we can decompose $u$ as

\begin{equation}  u |_{B_R} = \sum_\theta \sum_{T \in \TT_\theta} a_T W_T. \end{equation}

\noindent This is called the wave packet decomposition of $u$.

Wave packets originated in quantum physics.  Recall that when $S$ is a paraboloid, the image of $E_S$ is the set of solutions to the free Schrodinger equation, $\partial_t u = i \sum_{j=1}^{d-1} \partial_j^2 u$.  This free Schrodinger equation describes a quantum mechanical particle with no forces acting on it.  Now in some regimes, a quantum mechanical particle behaves like a classical particle.  Since there are no forces acting in our model, a classical particle would move at a constant velocity and trace out a straight line space-time.  A wave packet is a solution to the free Schrodinger equation that approximates a classical particle in this sense.  Physicists working on quantum mechanics noted that an arbitrary solution to the free Schrodinger equation can be written as a superposition of such wave packets.  This is the wave packet decomposition that Bourgain and others use in restriction theory.

The wave packet structure is useful in analysis because it constrains the geometry of the set where $| u_\theta |$ is large.  A basic fact about $W_T$ is that $|W_T|$ is roughly constant on $T$, and rapidly decaying away from $T$.  Therefore $|u_\theta(x)|$ is roughly constant on $x + \theta^*$, because $x+ \theta^*$ is essentially a single tube $T \in \TT_\theta$.   So if $|u_\theta(x)|$ is large, then $|u_\theta|$ must be large on the whole tube $x + \theta^*$.  Now if $|u(x)|$ is large, we might expect that $|u_\theta(x)|$ is large for many $\theta$.  Using orthogonality or the operator norm method, we can upper bound the volume of the set $\{ x \in B_R: |u_\theta(x)| > \lambda \}$.  Because of the wave packet structure, this set is a union of translates of $\theta^*$.  For different sets $\theta$, the tube $\theta^*$ points in different directions.  This makes it difficult for the sets $\{ x \in B_R : |u_\theta(x)| > \lambda \}$ to overlap too much, and so there are not many points $x$ where $|u_\theta(x)|$ is large for many $\theta$.  The estimates about how these sets overlap are called Kakeya type estimates.  In \cite{B3}, Bourgain used this strategy, together with the Tomas-Stein estimate, to prove new bounds for the restriction problem going beyond the $M M^*$ method and the power method.  

If $\tau \subset S$ is a cap of radius larger $r$ than $R^{-1/2}$, then $u_\tau = E_S (1_\tau g)$ also has wave packet structure.  In particular, $|u_\tau|$ is roughly constant on translates of $\tau^*$, a tube of radius $r^{-1}$ and length $r^{-2}$, with long direction perpendicular to the tangent space of $S$.  The proof of decoupling in \cite{BD} exploits this wave packet structure at all scales in a systematic way, together with orthogonality.  See \cite{Gu} for a survey article on the proof of decoupling.  


\newcommand{\image}{\textrm{Im}}

This wave packet decomposition is a special feature of the restriction problem and the matrix $M_S$.   For a random matrix $M_{ran}$, there is nothing remotely like wave packets.  Let us try to make this precise.  A cap $\theta$ corresponds to a subset of the columns of the matrix $M_S$.  If $M$ is any matrix and $\theta$ is a subset of the columns of $M$, then let us write $M_{;\theta}$ for the submatrix containing the columns in $\theta$.  (We use the semi-colon because we are writing $M_W$ for the submatrix containing the rows of $W$.  In general, we could write $M_{W;\theta}$ for the submatrix $M_{i,j}$ where $i \in W$ and $j \in \theta$.)   If we think of $M$ as a linear map $M: \CC^N \rightarrow \CC^T$, then the image $\image(M_{;\theta})$ is a subspace of $\CC^T$.  The wave packet structure is a very special feature of the subspace $\image(M_{S; \theta})$.  Suppose that $u_\theta = E_S (1_\theta g)$ as above -- morally $u_\theta$ is in $\image(M_{S; \theta})$.  Oversimplifying some technical details, we can describe the situation as follows.  Because of the wave packet structure of $u_\theta$ in (\ref{wptheta}), once we know the value of $u_\theta(x_0)$, we can determine the coefficent $a_T$ for the tube $T \in \TT_\theta$ containing $x_0$, and then we can determine $u_\theta(x)$ for all $x$ in this $T$.  So knowing $u_\theta(x_0)$ for a single $x_0$ essentially determines $u_\theta(x)$ for many other $x$.  In other words, for a vector $u_\theta$ in $\image(M_{S; \theta})$, knowing a single component of $u_\theta$ determines many other components.

Another way to say this is that there are very sparse vectors in the orthogonal complement of $\image(M_{S; \theta})$.  In our slightly oversimplified model above, there are vectors in this orthogonal complement with only two non-zero entries.  In the literal truth, there are vectors in this orthogonal complement with $O(1)$ entries of size $\sim 1$ and where the other entries have tiny norm.  This type of special structure does not occur for random matrices $M_{ran}$.  If $|\theta|$ denotes the cardinality of $\theta$, then $\image(M_{ran; \theta})$ usually has dimension $|\theta|$.  For any subspace of $\CC^T$ with dimension $|\theta|$, there is a vector in its orthogonal complement with at most $T - |\theta| +1$ non-zero entries, by simple linear algebra.  If the matrix $M$ is random, then I believe that for every $\theta$, the vectors in the orthogonal complement of $\image(M_{ran; \theta})$ are not much sparser than this.

For the Dirichlet polynomial matrix $M_{Dir}$, there are some wave packets, as discussed in Section \ref{secbourkak}, but the structure is not as useful as in the restriction problem.  The wave packet methods developed for the restriction problem starting with  \cite{B3} do not adapt to Dirichlet polynomials.  One way to see that is to consider the almost counterexample in Section \ref{subsecac}.  The wave packets of $M_{Dir}$ and of $M_{AC}$ are identical to each other.  So the input from wave packets cannot by itself distinguish the large value problem for Dirichlet polynomials from the almost counterexample.  In a large regime, our example involving $M_{AC}$ matches the bounds coming from the operator norm method and the $MM^*$ method, and so in this regime, the input from wave packets cannot by itself prove bounds going beyond those methods.  

(On a personal note, Yuqiu Fu, Dominique Maldague, and I spent a year studying Dirichlet polynomials using wave packets, leading to the paper \cite{FGM}.  We were able to use wave packets to prove strong bounds for short Dirichlet polynomials, of the form $\sum_{n = N}^{N + N^{1/2}} b_n n^{it}$, but our approach didn't give any new bounds for honest Dirichlet polynomials.  
Eventually, we found the almost counterexample from Section \ref{subsecac} and realized that our approach could not work for honest Dirichlet polynomials.)


This discussion suggests that sparse vectors in the orthogonal complement of $\image(M)$ and $\image( M_{;\theta})$ could be a useful tool to study large value estimates for $M$.  The wave packet methods we have described make use of certain sparse vectors in $\image(M_{S; \theta})^\perp$.  It would be interesting to know if there are any other sparse vectors in $\image(M_{S; \theta})^\perp$ or in $\image(M_{Dir; \theta})^\perp$ besides the ones coming from the wave packet structure.

\end{document}